\documentclass[a4paper]{article}
\pdfoutput=1

\usepackage{graphicx,wrapfig}
\usepackage{lipsum}
\usepackage{amsmath,amssymb,mathrsfs}
\usepackage{hyperref}
\usepackage{url}
\usepackage{mathtools}
\usepackage{changepage}
\usepackage{multirow}
\usepackage{wrapfig}
\usepackage{subfig}
\usepackage{color}
\usepackage{epsfig,enumerate,amsmath,amsfonts,amssymb,amsthm,mathrsfs,ifpdf}
\usepackage{indentfirst,relsize}
\usepackage{setspace,graphicx}
\usepackage{latexsym}
\usepackage[all]{xy}
\usepackage[usenames,dvipsnames]{pstricks}
\usepackage{pst-grad} 
\usepackage{pst-plot} 
\usepackage[margin = 2.50cm]{geometry}

\usepackage{algorithm}

\usepackage{algpseudocode}
\usepackage{color,soul}
\usepackage{authblk} 
\usepackage{enumitem}

\usepackage{orcidlink}

\newcommand{\floor}[1]{\left\lfloor #1 \right\rfloor}

\newcommand{\remove}[1]{}

\newtheorem{theorem}{Theorem}
\newtheorem{lemma}[theorem]{Lemma}
\newtheorem{proposition}[theorem]{Proposition}

\newtheorem{claim}[theorem]{Claim}
\newcounter{Case}[theorem]
\newtheorem{case}[Case]{Case}

\newtheorem{definition}[theorem]{Definition}
\newtheorem{remark}{Remark}
\newtheorem{observation}[theorem]{Observation}

\title{On total transitivity of graphs}

\author[]{Kamal Santra \orcidlink{0009-0006-5997-1452}\footnote{kamal.7.2013@gmail.com, kamal\_1821ma04@iitp.ac.in}}
\affil[]{Department of Mathematics\\
	
	Indian Institute of Technology Patna\\
	
	Bihta, 801106, Bihar, India}
\date{}

\begin{document}

\maketitle
\begin{abstract}
Let $G = (V, E)$ be a graph where $V$ and $E$ are the vertex and edge sets, respectively. For two disjoint subsets $A$ and $B$ of $V$, we say that $A$ \emph{dominates} $B$ if every vertex of $B$ is adjacent to at least one vertex of $A$. A vertex partition $\pi = \{V_1, V_2, \ldots, V_k\}$ of $G$ is called a \emph{transitive partition} of size $k$ if $V_i$ dominates $V_j$ for all $1 \leq i < j \leq k$. In this article, we study a variation of the transitive partition, namely the \emph{total transitive partition}. The total transitivity $Tr_t(G)$ is defined as the maximum order of a vertex partition $\pi = \{V_1, V_2, \ldots, V_k\}$ of $G$ obtained by repeatedly removing a total dominating set from $G$ until no vertices remain. Thus, $V_1$ is a total dominating set of $G$, $V_2$ is a total dominating set of the graph $G_1 = G - V_1$, and, for $2 \leq i \leq k - 1$, $V_{i+1}$ is a total dominating set in the graph $G_i = G - \bigcup_{j=1}^i V_j$. A vertex partition of order $Tr_t(G)$ is called a $Tr_t$-partition. The \textsc{Maximum Total Transitivity Problem} is to find a total transitive partition of a given graph with the maximum number of parts. First, we characterize split graphs with total transitivity equal to $1$ and $\omega(G) - 1$. Moreover, for a split graph $G$ and $1 \leq p \leq \omega(G) - 1$, we provide necessary conditions for $Tr_t(G) = p$. Furthermore, we show that the decision version of this problem is NP-complete for bipartite graphs. On the positive side, we prove that this problem can be solved in linear time for bipartite chain graphs. Finally, we design a polynomial-time algorithm to solve the \textsc{Maximum Total Transitivity Problem} in trees.
\end{abstract}


\section{Introduction}

Partitioning a graph is one of the fundamental problems in graph theory. In the partitioning problem, the objective is to partition the vertex set (or edge set) into some parts with desired properties, such as independence, minimal edges across partite sets, etc. In literature, partitioning the vertex set into certain parts so that the partite sets follow particular domination relations among themselves has been studied \cite{chang1994domatic, furedi2008inequalities, haynes2019transitivity, haynes2020upper, hedetniemi2018transitivity, paul2023transitivity, santra2023transitivity, zaker2005grundy}. Let $G$ be a graph with $V(G)$ as its vertex set and $E(G)$ as its edge set. When the context is clear, $V$ and $E$ are used instead of $V(G)$ and $E(G)$. The \emph{neighbourhood} of a vertex $v\in V$ in a graph $G=(V, E)$ is the set of all adjacent vertices of $v$ and is denoted by $N_G(v)$. The \emph{degree} of a vertex $v$ in $G$, denoted as $\deg_G(v)$, is the number of edges incident to $v$. A vertex $v$ is said to \emph{dominate} itself and all its neighbours. A \emph{dominating set} of $G=(V, E)$ is a subset of vertices $D$ such that every vertex $x\in V\setminus D$ has a neighbour $y\in D$, that is, $x$ is dominated by some vertex $y$ of $D$. 

There has been a lot of research on graph partitioning problems based on a domination relationship between the different sets. Cockayne and Hedetniemi introduced the concept of \emph{domatic partition} of a graph $G=(V, E)$ in 1977, in which the vertex set is partitioned into $k$ parts, say $\pi =\{V_1, V_2, \ldots, V_k\}$, such that each $V_i$ is a dominating set of $G$ \cite{cockayne1977towards}. The number representing the highest possible order of a domatic partition is called the \emph{domatic number} of G, denoted by $d(G)$. Another similar type of partitioning problem is \emph{Grundy partition}. Christen and Selkow introduced a Grundy partition of a graph $G=(V, E)$ in 1979 \cite{CHRISTEN197949}. In the Grundy partitioning problem, the vertex set is partitioned into $k$ parts, say $\pi =\{V_1, V_2, \ldots, V_k\}$, such that each $V_i$ is an independent set and for all $1\leq i< j\leq k$, $V_i$ dominates $V_j$. The maximum order of such a partition is called \emph{Grundy number} of $G$, denoted by $\Gamma(G)$. In 2018, Hedetniemi and Hedetniemi \cite{hedetniemi2018transitivity} introduced a transitive partition as a generalization of the Grundy partition. A \emph{transitive partition} of size $k$ is defined as a partition of the vertex set into $k$ parts, say $\pi =\{V_1,V_2, \ldots, V_k\}$, such that for all $1\leq i< j\leq k$, $V_i$ dominates $V_j$. The maximum order of such a transitive partition is called \emph{transitivity} of $G$ and is denoted by $Tr(G)$. In 2020, Haynes et al. generalized the idea of domatic partition as well as transitive partition and introduced the concept of \emph{upper domatic partition} of a graph $G$, where the set of the vertex is divided into $k$ parts, say $\pi =\{V_1, V_2, \ldots, V_k\}$, such that for each $i, j$, with $1\leq i<j\leq k$, either $V_i$ dominates $V_j$ or $V_j$ dominates $V_i$ or both \cite{haynes2020upper}. The maximum order of such an upper domatic partition is called \emph{upper domatic number} of $G$, denoted by $D(G)$. All these problems, domatic number \cite{chang1994domatic, zelinka1980domatically, zelinka1983k}, Grundy number \cite{effantin2017note, furedi2008inequalities, hedetniemi1982linear, zaker2005grundy, zaker2006results}, transitivity \cite{haynes2019transitivity, hedetniemi2018transitivity, paul2023transitivity, santra2023transitivity}, upper domatic number \cite{haynes2020upper, samuel2020new} have been extensively studied both from an algorithmic and structural point of view. A Grundy partition is a transitive partition with the additional restriction that each partite set must be independent. In a transitive partition $\pi =\{V_1, V_2, \ldots, V_k\}$ of $G$, we require the domination property in one direction, that is, $V_i$ dominates $V_j$ for $1\leq i< j\leq k$. However, in an upper domatic partition $\pi =\{V_1,V_2, \ldots, V_k\}$ of $G$, for all $1\leq i<j\leq k$, either $V_i$ dominates $V_j$ or $V_j$ dominates $V_i$ or both. The definition of each vertex partitioning problem ensures the following inequalities for any graph $G$. For any graph $G$, $1\leq \Gamma(G)\leq Tr(G)\leq D(G)\leq n$, where $n$ is the number of vertices of $G$.

\subsection{Motivation}
Many variations of the classical domination problem have been studied in the literature \cite{haynes2020topics, haynes2021structures, haynes2013fundamentals}. One natural variation is the \emph{total dominating set}, where the goal is to find a set $D$ such that every vertex in $V$ is dominated by a vertex in $D$. The concept of total domination was first formally introduced in 1980 by Cockayne et al. \cite{cockayne1980total} as part of their study on domination in graphs.In the classical domination problem, the objective is to monitor the vertices of a network that are not in a dominating set by using the vertices of the dominating set. Total domination ensures that every node in a network is adjacent to a control node, making it useful for designing fault-tolerant systems where every node has a backup connection.This variation of the domination problem has been well-studied in the literature from both structural and algorithmic perspectives. It has been shown that finding the minimum total dominating set is NP-complete for general graphs, bipartite graphs \cite{pfaff1983np}, and split graphs \cite{laskar1983domination}, whereas it is solvable in polynomial time for trees \cite{laskar1984algorithmic}.

The domatic number of a graph was introduced by Cockayne and Hedetniemi \cite{cockayne1977towards}. Furthermore, many variants of this concept were introduced and studied by B. Zelinka \cite{zelinka1992domatic}. One such variation is the total domatic number, which was introduced by Cockayne et al. \cite{cockayne1980total}. Several works have been done on the total domatic number \cite{zelinka1988total, zelinka1989total}. In 2024, Paul and Santra introduced three variations of transitivity, namely, $2$-transitivity \cite{PAUL202457}, $d_2$-transitivity \cite{santra2024d_2transitivity}, and strong transitivity \cite{santra2024strong}. The $2$-transitivity, $d_2$-transitivity, and strong transitivity are based on $2$-domination, $d_2$-domination, and strong domination, respectively. In this article, we explore another natural variation of the transitivity problem, focusing on total domination instead of domination. We introduce such a variation of the transitivity problem, namely \emph{total transitive partition}. Total transitivity $Tr_t(G)$ equals the maximum order of a vertex partition $\pi = \{V_1, V_2, \ldots, V_k\}$ of $G$ obtained by repeatedly removing a total dominating set from $G$ until no vertices remain. Thus, $V_1$ is a total dominating set of $G$, $V_2$ is a total dominating set of the graph $G_1=G-V_1$, and, for $2\leq i\leq k-1$, $V_{i+1}$ is a total dominating set in the graph $G_i=G-\bigcup_{j=1}^i V_i$. A vertex partition of order $Tr_t(G)$ is called a $Tr_t$-partition. Note that every total transitive partition is also a transitive partition. Therefore, for any graph $G$, $1\leq Tr_t(G)\leq Tr(G)\leq n$. Since the total dominating set exists only when $G$ does not have any isolated vertices, in this paper, we focus only on graphs with no isolated vertices.

An equivalent definition of total transitivity can be given in terms of total domination among partite sets and the set itself. A \emph{total transitive partition} of size $k$ is defined as a partition of the vertex set into $k$ parts, say $\pi =\{V_1,V_2, \ldots, V_k\}$, such that for all $1\leq i\leq j\leq k$, $V_i$  totally dominates $V_j$. The maximum order of such a total transitive partition is called the \emph{total transitivity} of $G$ and is denoted by $Tr_t(G)$. 

We define another variation of transitivity, namely \emph{weak total transitivity}. A \emph{weak total transitive partition} of size $k$ is defined as a partition of the vertex set into $k$ parts, say $\pi =\{V_1,V_2, \ldots, V_k\}$, such that for all $1\leq i\leq j\leq k$ and $i\neq k$, $V_i$ totally dominates $V_j$. The maximum order of such a weak total transitive partition is called the \emph{weak total transitivity} of $G$ and is denoted by $Tr_t^w(G)$. Later, we will use this parameter to solve total transitivity in a tree. Clearly, if $\pi$ is a total transitive partition, it is also a weak total transitive partition, and hence, for a graph $G$, $Tr_t(G)\leq Tr_t^w(G)$. Also, note that in a weak total transitive partition $\pi=\{V_1, V_2, \ldots, V_k\}$, $V_k$ may not be a total dominating set of $G-\bigcup_{j=1}^{k-1}V_i$  but if $\pi$ is a total transitive partition, $V_k$ is a total dominating set of $G-\bigcup_{j=1}^{k-1}V_i$.

The \textsc{Maximum Total Transitivity Problem} and its corresponding decision version are defined as follows.

\begin{center}

	\fbox{%
		\parbox{0.7\linewidth}{%
			\noindent\textsc{Maximum Total Transitivity Problem}
			
			\noindent\emph{Instance:} A graph $G=(V,E)$
			
			\noindent\emph{Solution:} A total transitive partition of $G$ with maximum size%
		}%
	}
	
	\vspace{0.5cm}

	\fbox{%
		\parbox{0.8\linewidth}{%
			\noindent\textsc{Maximum Total Transitivity Decision Problem}
			
			\noindent\emph{Instance:} A graph $G=(V,E)$, integer $k$
			
			\noindent\emph{Question:} Does $G$ have a total transitive partition of order at least $k$?%
		}%
	}
	
\end{center}

In this paper, we study this variation of transitive partition from a structure and algorithmic point of view. The main contributions are summarized below:

\begin{enumerate}[label=(\arabic*.)]
	
	\item We characterize split graphs with total transitivity equal to $1$ and $\omega(G)-1$.
	
	\item Given a split graph $G$ and $1 \leq p \leq \omega(G) - 1$, we provide necessary conditions for $Tr_t(G) = p$.
	
	\item We show that the \textsc{Maximum Total Transitivity Decision Problem} is NP-complete for bipartite graphs.
	
	\item We prove that \textsc{Maximum Total Transitivity Problem} can be solved in linear time for bipartite chain graphs.
	
	\item We design a polynomial time algorithm for \textsc{Maximum Total Transitivity Problem} in trees.

\end{enumerate}

The rest of the paper is organized as follows. Section 2 contains basic definitions and notations that are followed throughout the article. This section also discusses the properties of total transitivity of graphs. In section 3, we characterize split graphs with total transitivity equal to $1$ and $\omega(G)-1$. Also, in this section, we provide necessary conditions for $Tr_t(G) = p$ for split graphs $G$ and $1\leq p\leq \omega(G)-1$. Section 4 shows that the \textsc{Maximum Total Transitivity Decision Problem} is NP-complete in bipartite graphs. In section 5 we prove that \textsc{Maximum Total Transitivity Problem} can be solved in linear time for bipartite chain graphs. In Section 6, we design a polynomial-time algorithm for solving \textsc{Maximum Total Transitivity Problem} in trees.  Finally, Section 7 concludes the article with some open problems.

\section{Preliminaries}

\subsection{Notation and definition}

Let $G=(V, E)$ be a graph with $V$ and $E$ as its vertex and edge sets, respectively. A graph $H=(V', E')$ is said to be a \emph{subgraph} of a graph $G=(V, E)$ if and only if $V'\subseteq V$ and $E'\subseteq E$. For a subset $S\subseteq V$, the \emph{induced subgraph} on $S$ of $G$ is defined as the subgraph of $G$ whose vertex set is $S$ and whose edge set consists of all of the edges in $E$ that have both endpoints in $S$, and it is denoted by $G[S]$. The \emph{complement} of a graph $G=(V,E)$ is the graph $\overline{G}=(\overline{V}, \overline{E})$, such that $\overline{V}=V$ and $\overline{E}=\{uv| uv\notin E \text{ and } u\neq v\}$. The \emph{open neighbourhood} of a vertex $x\in V$ is the set of vertices $y$ adjacent to $x$, denoted by $N_G(x)$. The \emph{closed neighborhood} of a vertex $x\in V$, denoted as $N_G[x]$, is defined by $N_G[x]=N_G(x)\cup \{x\}$. A subset $D\subseteq V$ is said to be a \emph{total dominating set} of $G$, if for every vertex $x\in V$, there exists a vertex $y\in D$ such that $xy\in E$. \emph{Breadth first search} (BFS) is a graph searching technique that systematically explores all the vertices of a connected graph. The process starts at a designated vertex and explores all its neighbours before moving on to the next level of neighbours. The order in which the vertices are visited is called the \emph{BFS ordering}. A vertex is said to be a \emph{leaf vertex} if its degree is $1$, and the vertex that is adjacent to a leaf vertex is said to be the \emph{support vertex} of the leaf vertex.

A subset of $S\subseteq V$ is said to be an \emph{independent set} of $G$ if no two vertices in $S$ are adjacent. A subset of $K\subseteq V$ is said to be a \emph{clique} of $G$ if every pair of vertices in $K$ is adjacent. The cardinality of a maximum size clique is called \emph{clique number} of $G$, denoted by $\omega(G)$.  A graph is called \emph{bipartite} if its vertex set can be partitioned into two independent sets. A graph $G=(V, E)$ is a \emph{split graph} if $V$ can be partitioned into an independent set $S$ and a clique $K$.

A bipartite graph $G=(X\cup Y,E)$ is called a \textit{bipartite chain graph} if there exists an ordering of vertices of $X$ and $Y$, say $\sigma_X= (x_1,x_2, \ldots ,x_{n_1})$ and $\sigma_Y=(y_1,y_2, \ldots ,y_{n_2})$, such that $N(x_{n_1})\subseteq N(x_{n_1-1})\subseteq \ldots \subseteq N(x_2)\subseteq N(x_1)$ and $N(y_{n_2})\subseteq N(y_{n_2-1})\subseteq \ldots \subseteq N(y_2)\subseteq N(y_1)$. Such ordering of $X$ and $Y$ is called a \emph{chain ordering}, and it can be computed in linear time \cite{heggernes2007linear}.

\subsection{Basic properties of total transitivity}
In this subsection, we present some properties of total transitivity that also motivate us to study this variation of transitivity. First, we show the following bound for total transitivity.

\begin{proposition} \label{upper_bound_total_transitivity}
	For any graph $G$, $Tr_t(G)\leq \min\{\Delta(G), \floor{\frac{n}{2}}\}$, where $\Delta(G)$ is the maximum degree of $G$.
\end{proposition}

\begin{proof}
	Let $\pi=\{V_1, V_2, \ldots, V_k\}$ be a $Tr_t$-transitive partition of $G$. Let $x\in V_k$. As $V_1$ is a total dominating set of $G$, and for $1 \leq i \leq k - 1$, $V_{i+1}$ is a total dominating set in the graph $G_i = G - \bigcup_{j=1}^i V_j$, $deg(x)\geq k$ and $V_j\geq 2$ for all $1\leq j\leq k$. Therefore, $Tr_t(G)=k\leq deg(x)\leq \Delta(G)$ and $Tr_t(G)\leq \floor{\frac{n}{2}}$. Hence, $Tr_t(G)\leq \min\{\Delta(G), \floor{\frac{n}{2}}\}$.
\end{proof}

Next, we find the total transitivity of complete graphs, paths, cycles and complete bipartite graphs.
\begin{proposition} \label{Complete_graph_total_transitivity}
	For the complete graph $K_n$ of $n$ vertices $Tr_t(K_n)=\floor{\frac{n}{2}}$.
\end{proposition}
\begin{proof}
	The partition $\{V_1, V_2, \ldots, V_{\floor{\frac{n}{2}}}\}$ such that each set contains at least two vertices is a total transitive partition of $K_n$ of size $\floor{\frac{n}{2}}$. Moreover, from the Proposition \ref{upper_bound_total_transitivity}, we have $Tr_t(G)\leq \min\{\Delta(K_n), \floor{\frac{n}{2}}\}= \min\{n,\floor{\frac{n}{2}}\}=\floor{\frac{n}{2}}$. Therefore, $Tr_t(K_n)=\floor{\frac{n}{2}}$.
\end{proof}

\begin{remark}\label{rem_total_arbitrarily}
	From the definition of a total transitive partition, it is clear that it is also a transitive partition of $G$. But the difference $Tr(G)-Tr_t(G)$ can be arbitrarily large. From the Proposition \ref{Complete_graph_total_transitivity}, we have $Tr_t(K_n)=\frac{n}{2}$ for $n$ is an even integer, and from \cite{hedetniemi2018transitivity}, we know that $Tr(K_n)=n$, hence the difference $Tr(K_n)-Tr_t(K_n)=\frac{n}{2}$ is arbitrarily large for $n$. 
\end{remark}

\begin{proposition}\label{Path_total_tr}
	Let $P_n$ be a path of $n$ vertices, and then the total transitivity of $P_n$ is given as follows:
	
	$Tr_t(P_n) = \begin{cases}
		1 & n=2, 3, 4, 5 \\
		2 & n\geq 6
	\end{cases}$
\end{proposition}

\begin{proposition}\label{Cycle_total_tr}
	Let $C_n$ be a cycle of $n$ vertices, then the total transitivity of $C_n$ is given as follows:
	
	$Tr_t(C_n) = \begin{cases}
		1 & n=3 \\
		2 & n\geq 4
	\end{cases}$
\end{proposition}

\begin{proposition} \label{Complete_bipartite_graph_ttr}
	For a complete bipartite graph $K_{m,n}$ with $m+n$ vertices $Tr_t(K_{m,n})=\min\{m, n\}$.

\end{proposition}
\begin{proof}
	Let the vertex set of $K_{m,n}$ be $V = X \cup Y$, where $|X| = m$, $|Y| = n$, and assume $m \leq n$. Consider the vertex partition $\pi = \{V_1, V_2, \dots, V_m\}$, where $V_i = \{x_i, y_i\}$ for all $2 \leq i \leq m$, and $V_1$ contains all remaining vertices of $K_{m,n}$. Since $K_{m,n}$ is a complete bipartite graph, this partition is a total transitive partition of size $m$, implying that $Tr_t(K_{m,n}) \geq m$.
	
	Next, we show that $Tr_t(K_{m,n})$ cannot exceed $m$. Suppose, for contradiction, that $Tr_t(K_{m,n}) = k \geq m+1$, and let $\pi' = \{V_1, V_2, \dots, V_k\}$ be a total transitive partition of $K_{m,n}$ with size $k$. Choose any vertex $v \in V_k$. Since $v$ must be totally dominated within $V_k$, it follows that the degree of $v$ must be at least $k \geq m+1$. However, this condition is only possible if $v \in X$, as vertices in $Y$, have a degree at most $m$. Moreover, since $V_k$ totally dominates itself, there must exist $u \in V_k$ such that $uv \in E$ and $u \in Y$. Since $u \in V_k$, it must also have a degree at least $k \geq m+1$. But this contradicts the fact that $\deg(u) = m < k$. Hence, no such partition of size greater than $m$ can exist, implying that $Tr_t(K_{m,n}) \leq m$. Thus, we conclude that $Tr_t(K_{m,n}) =\min\{m, n\}$.
\end{proof}

Similarly, as transitivity, if we union any two sets of a total transitive partition of order $k$, we create a total transitive partition of order $k-1$. Hence, we have a similar proposition as transitivity.

\begin{proposition}
	If $\pi=\{V_1, V_2, \ldots, V_k\}$ is a total transitive partition of a graph $G$, then for any two sets $V_i$ and $V_j$ in $\pi$, where $i<j$, the partition $\pi'=\{V_1, V_2, \ldots, V_{i-1}, V_i\cup V_j, V_{i+1}, \ldots, V_{j-1}, V_{j+1}, \ldots, V_k\}$ is a total transitive partition.  
\end{proposition}
\begin{proof}
	By the definition of total transitive partition, clearly it follows that $\pi'$ is a total transitive partition. of $G$ of size $k-1$.
\end{proof}

Given this, we have the following interpolation result.

\begin{proposition}\label{any_size_totaltr}
	Let $G$ be a graph and $k$ be the order of a total transitive partition of $G$, then for every $j$, $1\leq j\leq k$, $G$ has a total transitive partition of order $j$.
\end{proposition}


For disconnected graph we have the following proposition.

\begin{proposition}\label{disjoint_union_total_transitivity}
	Let $G$ be a disconnected graph with connected components $C_1, C_2, \ldots, C_t$. Then $Tr_t(G)=\max\{Tr_t(C_i) | 1\leq i\leq k\}$.
\end{proposition}

\begin{proof}
	Let $C_j$ be the component of $G$ such that $\max\{Tr_t(C_i) | 1\leq i\leq k\}=Tr_t(C_j)=p$. Furthermore, assume that $\pi=\{V_1, V_2, \ldots, V_p\}$ is a total transitive partition of $C_j$ of size $p$. Consider $\pi'=\{(\bigcup_{r=1, r\neq j}^t C_r)\cup V_1, V_2, \ldots, V_p\}$. Clearly, $\pi'$ is a total transitive partition of $G$ of size $p$ and hence $Tr_t(G)\geq p=Tr_t(C_j)=\max\{Tr_t(C_i) | 1\leq i\leq k\}$.
	
	For the other part, let $\pi=\{V_1, V_2, \ldots, V_k\}$ be a $Tr_t$-partition of $G$. Furthermore, let for some $l\in\{1, 2, \ldots, k\}$, $V(C_l)\cap V_k\not=\phi$. Since $V_j$ totally dominates $V_k$ for all $1\leq j\leq k$, then $V(C_l)\cap V_j\not=\phi$. For a similar reason, we have $(V_j \cap V(C_l))$, which totally dominates $(V_k \cap V(C_l))$. Therefore, $\{V(C_l)\cap V_1, V(C_l)\cap V_2, \ldots, V(C_l)\cap V_k \}$ is a total transitive partition of $C_l$. Thus, $k=Tr_t(G)\leq Tr_t(C_l)$, that is, $Tr_t(G)\leq \max\{Tr_t(C_i) | 1\leq i\leq k\}$. Hence, $Tr_t(G)= \max\{Tr_t(C_i) | 1\leq i\leq k\}$.
\end{proof}

Now we show that if the total transitivity of a graph $G$ is constant, that is, independent of the number of vertices, and $\pi=\{V_1, V_2, \ldots, V_{k-1}, V_k\}$ is a total transitive partition of $G$, then the number of vertices of each $V_i$, for $2\leq i\leq k$, does not depend on $n$; rather, it is bounded by some exponential in $i$.
\begin{proposition}\label{size of total_transitive_partition}
	Let $G$ be a connected graph and $Tr_t(G)=k$, $k\geq 3$. Then there exists a total transitive partition of $G$ of size $k$, say $\pi =\{V_1, V_2, \ldots, V_k\}$ of $G$, such that $|V_k|=2$ and $|V_{k-i}|\leq 2^{i+1}$, for all $i$, $1\leq i\leq k-2$. 
\end{proposition}
\begin{proof}
	Let $\pi =\{V_1, V_2, \ldots, V_k\}$ be a total transitive partition of $G$ of size $k$, where $k\geq 3$. Let $x\in V_k$; as $V_k$ totally dominates $V_k$, there exists $y\in V_k$ such that $xy\in E(G)$. So, $|V_k|\geq 2$. Now if $|V_k-\{x, y\}|\geq 1$, then we move all the vertices except two, say vertices $x$ and $y$, from the set $V_k$ to the first set $V_1$. After this transformation, we have $\pi' =\{V_1'=V_1\cup (V_k-\{x, y\}), V_2, \ldots, V_{k-1}, V_k'=\{x, y\}\}$, which is a total transitive partition of $G$ of size $k$ and $|V_k|=2$. Therefore, $\pi'$ is a total transitive of $G$ of size $k$, such that $|V_k|=2$. Consider the base case when $i=1$. Note that to totally dominate $V_k$ and $V_{k-1}$, we need at most $2\cdot|V_k|$ vertices in $V_{k-1}$, and the remaining vertices can be moved to $V_1$ as in the previous argument. Therefore, $|V_{k-1}|\leq 2\cdot|V_k|=2\cdot 2=4=2^{i+1}$. Assume the induction hypothesis and consider that we have a total transitive partition $\pi=\{V_1, V_2, \ldots\, V_k\}$, such that $|V_k|=2$ and $|V_{k-j}|\leq 2^{j+1}$, $1\leq j\leq i-1$. Now consider the set $V_{k-i}$. Similarly, at most $2\cdot |V_{k-i+1}|$ vertices in $V_{k-i}$ are needed to totally dominate the $i+1$ sets that follow and include it. If $V_{k-i}$ contains more than that many vertices, the remaining vertices can be moved to $V_1$, maintaining the total transitivity of $G$. Hence, $V_{k-i}\leq 2\cdot|V_{k-i+1}|= 2\cdot 2^{i}=2^{i+1}$.
\end{proof}

%
%
%
%
%
%

\begin{figure}[htbp!]
	\centering
	\includegraphics[scale=0.80]{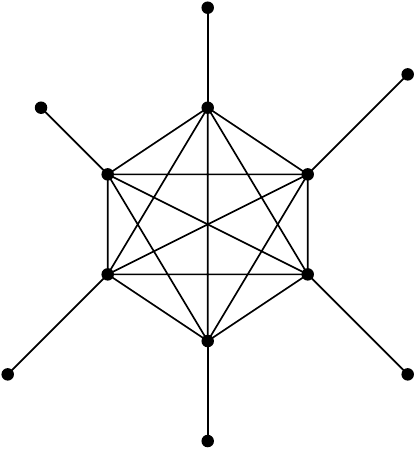}
	\caption{Example of a split graph where $Tr_t(G)=1$ and $Tr(G)=\omega(G)+1=6$}
	\label{fig:split_graph_example}
\end{figure}

In the Remark \ref{rem_total_arbitrarily}, we have already shown that the difference between transitivity and total transitivity can be arbitrarily large. Now, we show another non-trivial example where the difference between these two parameters is arbitrarily large. Let $G=(S\cup K, E)$ be a split graph, where $S$ and $K$ are an independent set and clique of $G$, respectively. Furthermore, assume that $K$ is the maximum clique of $G$, that is, $\omega(G)=|K|$. Furthermore, $|K|=|S|$, and each vertex of $S$ is adjacent to exactly one vertex of $K$ (see Figure \ref{fig:split_graph_example}).

First, we show that for the split graph given in Figure \ref{fig:split_graph_example}, $Tr_t(G)=1$. Let $\pi=\{V_1, V_2, \ldots, V_k\}$ be a total transitive partition of the above split graph $G$. Since $deg(s)=1$ for all $s\in S$, $s$ cannot be in $V_p$ for $p\geq 2$. Moreover, $\pi$ is a total transitive partition, which implies that to totally dominate $s\in S$, there must be a vertex from $K$ in $V_1$. From the above example graph, it is clear that to totally dominate vertices from $S$, all the vertices of $K$ must be in $V_1$. Hence, $Tr_t(G)=1$. Again, we know from \cite{santra2023transitivity}, $Tr(G)=\omega(G)+1$, as $S$ dominates $K$. Therefore, for large $\omega(G)$, the difference $Tr(G)-Tr_t(G)=\omega(G)+1-1=\omega(G)$ is arbitrarily large.

Next, we prove the existence of graphs with transitivity is one more than total transitivity. We define a special tree, namely \emph{total transitive tree of order $k$} and we denote it $\mathcal{T}_k$. Apart from this property, the $\mathcal{T}_k$ possesses some important properties. The recursive definition of $\mathcal{T}_k$ is as follows:

\begin{definition}
	The total transitive tree of order $k$, denoted as $\mathcal{T}_k$, is defined recursively as follows:
	\begin{enumerate}
		\item The total transitive tree of order $1$, denoted as $\mathcal{T}_1$, is a rooted tree, whose root has degree $1$ and adjacent to a vertex of degree $1$.
		
		
		\item The total transitive tree of order $k$, denoted as $\mathcal{T}_k$, is a rooted tree whose root has degree $k$ and is adjacent to vertices $\{v_1, v_2, \ldots, v_k\}$. For all $1\leq i \leq k-1$, $v_i$ are roots of the total transitive tree of order $i$, denoted by $\mathcal{T}_i$. Moreover, the vertex $v_k$ is adjacent to vertices $\{u_1, u_2, \ldots, u_{k-1}\}$, such that $u_i$ are roots of the total transitive tree of order $i$ ($\mathcal{T}_i$), for all $1\leq i \leq k-1$.
	\end{enumerate}
\end{definition}

\begin{figure}[htbp!]
	\centering
	\includegraphics[scale=0.80]{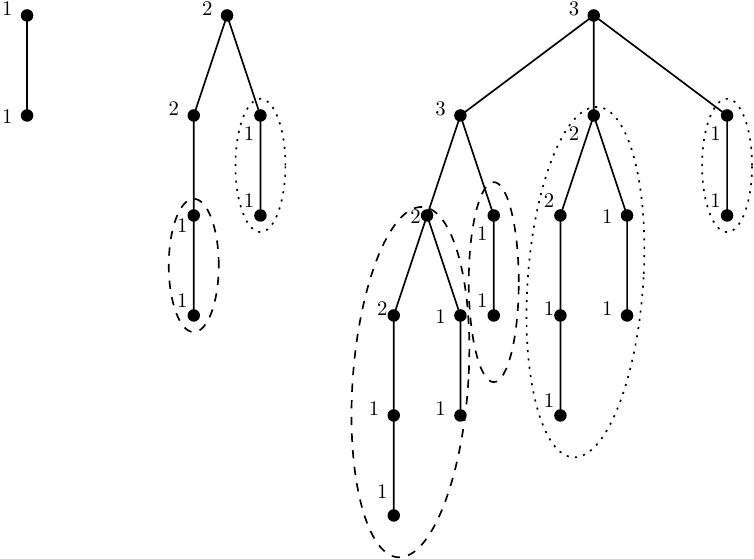}
	\caption{Total transitive tree of order up to $3$ along with its unique total transitive partition}
	\label{fig:first_three_tcmbt}
\end{figure}

Figure~\ref{fig:first_three_tcmbt} illustrates the total transitive tree of order up to $3$. Note that from the definition, it follows that the root of $\mathcal{T}_k$ and exactly one of its children are the only vertices that can be in $V_k$, in a total transitive partition of order $k$. Because of the acyclic property, none of the vertices in $V_i$ can dominate more than one vertex from $V_j$ with $j>i$. So, $\mathcal{T}_k$ is the smallest tree (in terms of vertices) with total transitivity equal to $k$. It also has the only total transitive partition of size $k$ (see Figure \ref{fig:first_three_tcmbt}). This property of $\mathcal{T}_k$ helps us to use these graphs as a gadget in the NP-completeness proofs in Section \ref{Section_3_total}. We show that if a graph $G$ is a $\mathcal{T}_k$, then $Tr_t(G)=k$ and $Tr(G)=k+1$. For that, we have the following lemma.
\begin{lemma}\label{total_transitivity_of_tcmbt_k}
	Let $T$ be a $\mathcal{T}_k$. Then $Tr_t(T)=k$ and $Tr(T)=k+1$.
\end{lemma}

\begin{proof}
	First, we show that $Tr_t(T)=k$. For this, we use induction on the order of the total transitive tree. When $k=1$, the base case, $T$, is a $P_2$ by the definition of a total transitive tree. Therefore, $Tr_t(P_2)=1$. Let $T$ be a $\mathcal{T}_k$. By the definition, we have $T$ is a rooted tree, whose root has degree $k$ and is adjacent to vertices $\{v_1, v_2, \ldots, v_k\}$, such that $v_i$ are roots of the total transitive tree of order $i$ $(\mathcal{T}_i)$, for all $1\leq i \leq k-1$. Moreover, the vertex $v_k$ is adjacent to vertices $\{u_1, u_2, \ldots, u_{k-1}\}$, such that $u_i$ are roots of the total transitive tree of order $i$ $(\mathcal{T}_i)$, for all $1\leq i \leq k-1$. Let $c$ be the root of $T$ and $c'$ be its child such that $c'$ is adjacent to vertices $\{u_1, u_2, \ldots, u_{k-1}\}$, such that $u_i$ are roots of the total transitive tree of order $i$ $(\mathcal{T}_i)$, for all $1\leq i \leq k-1$. By the induction hypothesis, $Tr_t(T_{v_{j}})=Tr_t(T_{u_{j}})=j$, for all $1\leq i \leq k-1$, where $T_v$ is a subtree of $T-\{c\}$, rooted at $v$ (see the $\mathcal{T}_3$ in Figure \ref{fig:first_three_tcmbt}). Now $T-\{c, c'\}$ is a disjoint union of trees $T_{v_j}, T_{u_j} $, for all $1\leq j\leq k-1$, and $k-1=\max\{Tr_t(T_{v_j})|1\leq j\leq k-1\}=\max\{Tr_t(T_{u_j})|1\leq j\leq k-1\}$. Total transitivity depends only on the maximum total transitivity of a component of a disconnected graph (Proposition \ref{disjoint_union_total_transitivity}). Therefore, $Tr_t(T-\{c, c'\})$ is $k-1$. Moreover, we can construct a total transitive partition, say $\pi=\{V_1, V_2, \ldots, V_{k-1}\}$, of $T-\{c, c'\}$ recursively such that $\{v_i, u_i\}\in V_i$ for all $1\leq i \leq k-1$. Now clearly, $\pi'= \pi \cup \{c, c'\}$ forms a total transitive partition of $T$. Therefore, $Tr_t(T)\geq k$. Moreover, by the Proposition \ref{upper_bound_total_transitivity}, we know that $Tr_t(T)\leq \Delta(G)=k$. Hence, $Tr_t(T)=k$. The algorithm described in \cite{hedetniemi1982linear} shows that $Tr(T)=k+1$. Therefore, $Tr_t(T)=k$ and $Tr(T)=k+1$ for $\mathcal{T}_k$.
\end{proof}

\section{Total Transitivity in split graphs}

In this section, we show that total transitivity of a split graph can be between $1$ and $\omega(G)-1$, where $\omega(G)$ is the size of a maximum clique in $G$. Then, we characterize split graphs with total transitivity of $1$ and $\omega(G)-1$. Later, we show some necessary conditions for $Tr_t(G)=p$, $1\leq p\leq \omega(G)-1$. First, we prove that total transitivity of a split graph $G$ cannot be greater than $\omega(G)-1$.

\begin{lemma}\label{split_graph_total_transitivity_lemma_1}
	Let $G=(S\cup K, E)$ be a split graph, where $S$ and $K$ are an independent set and clique of $G$, respectively. Also, assume that $K$ is the maximum clique of $G$, that is, $\omega(G)=|K|$. Then $1\leq Tr_t(G)\leq \omega(G)-1$.
\end{lemma}
\begin{proof}
	Clearly, $Tr_t(G)\geq 1$. Also, from the definition of total transitivity, we know that $Tr_t(G)\leq Tr(G)$. Again, from \cite{santra2023transitivity}, we know $Tr(G)\leq \omega(G)+1$. Therefore, for a split graph, $Tr_t(G)\leq \omega(G)+1$. Now we show that $Tr_t(G)$ cannot be $\omega(G)+1$ or $\omega(G)$. Let $Tr_t(G)=\omega(G)+1$ and $\pi=\{V_1, V_2, \ldots, V_{\omega(G)+1}\}$ be a total transitive partition of $G$. Since $|K|=\omega(G)$, there exists at least one set in $\pi$, say $V_i$, such that $V_i$ contains only vertices from $S$. Note that, in this case, $V_i$ cannot totally dominate $V_i$ as $S$ is an independent set of $G$. Therefore, we have a contradiction. Hence, $Tr_t(G)<\omega(G)+1$.
	
	Further assume $Tr_t(G)=\omega(G)$ and $\pi=\{V_1, V_2, \ldots, V_{\omega(G)}\}$ is a total transitive partition of $G$. By the Proposition \ref{size of total_transitive_partition}, we can assume that $|V_{\omega(G)}|=2$. Now $V_{\omega(G)}$ contains either two vertices from $K$ or one vertex from $K$ and another vertex from $S$. For the first case, there exists at least one set in $\pi$, say $V_i$, such that $V_i$ contains only vertices from $S$. Because $S$ is independent, $V_i$ cannot totally dominate $V_i$, so we have a contradiction. Consider the other case, namely that $V_{\omega(G)}$ contains one vertex from $K$ and another vertex from $S$. Let $s\in S\cap V_{\omega(G)}$. Then, by the definition of total transitive partition, $deg(s) \geq \omega(G)$. But $|K|=\omega(G)$ implies that $deg(s) < \omega(G)$, which is a contradiction. Therefore, in each case, we have a contradiction. Hence, $1\leq Tr_t(G)\leq \omega(G)-1$.
\end{proof}

Next, we characterize split graphs with $Tr_t(G)=1$. A vertex is said to be a \emph{leaf vertex} if its degree is $1$, and the vertex that is adjacent to a leaf vertex is said to be the \emph{support vertex} of the leaf vertex.

\begin{theorem}\label{split_graph_total_transitivity_characterization_equals_one}
	Let $G=(S\cup K, E)$ be a split graph, where $S$ and $K$ are the independent set and clique of $G$, respectively. Also, assume that $K$ is the maximum clique of $G$, that is, $\omega(G)=|K|$. Then $Tr_t(G)=1$ if and only if $G$ satisfies one of the following conditions.
	
	\begin{enumerate}[label=(\alph*)]
		\item For all $s\in S$, $deg(s)=1$, and $K$ has at most one vertex without a neighbour in $S$.
		
		\item For all $s\in S$ with $deg(s)\geq 2$, each vertex from $N(s)$ is a support vertex, and $K$ has at most one vertex without a neighbour in $S$.
	\end{enumerate}
\end{theorem}

\begin{proof}
	
	Let G satisfy the condition (a). We show that $Tr_t(G)=1$. Assume $\pi=\{V_1, V_2, \ldots, V_t\}$ is a total transitive partition of $G$. Since $deg(s)=1$ for all $s\in S$, no vertices from $S$ can be in $V_p$ for $p\geq 2$. To totally dominate vertices from $S$, $V_1$ contains all vertices from $K$ except at most one. Thus $V_2$ can contain at most one vertex from $K$, which contradicts the fact that $\pi$ is a total transitive partition. Therefore, $t$ cannot be more than $2$. Hence, $Tr_t(G)=1$.
	
	For condition $(b)$, let $\pi = \{V_1, V_2, \dots, V_t\}$ be a total transitive partition of $G$. Clearly, if $deg(s) = 1$ for some $s \in S$, then $s$ must belong to $V_1$. Now, assume there exists $s' \in S$ such that $deg(s') \geq 2$. By condition $(b)$, each vertex in $N(s')$ is a support vertex. Therefore, each vertex in $N(s')$ must be in $V_1$, and hence $s'$ cannot be in $V_p$ for $p \geq 2$. Thus, every vertex in $S$ must be in $V_1$, and at most one vertex from $K$ can be in $V_1$. This contradicts the fact that $\pi$ is a total transitive partition. Therefore, $t$ cannot be greater than 2. Hence, $Tr_t(G) = 1$.

	Conversely, assume that $Tr_t(G)=1$. We divide our proof into two cases:
	
	\begin{case}\label{Lower_STTCA1}
		The $deg(s)=1$ for all $s\in S$
	\end{case}
	We show that $K$ has at most one vertex without a neighbour in $S$. If $K$ has more than one vertex without a neighbour in $S$, we will derive a contradiction. Let $\{k_1, k_2\} \subseteq K$ such that $k_1$ and $k_2$ have no neighbour in $S$. Consider a vertex partition $\{V_1, V_2\}$ of $G$ as follows: $V_2 = \{k_1, k_2\}$, and $V_1$ contains all the other vertices of $G$. Clearly, $\{V_1, V_2\}$ is a total transitive partition of $G$ with size 2, which implies that $Tr_t(G) \geq 2$, a contradiction since $Tr_t(G) = 1$. Hence, $K$ has at most one vertex without a neighbour in $S$, and $G$ satisfies condition $(a)$.

	\begin{case}\label{Lower_STTCA2}
		There exists $s\in S$ with $deg(s)\geq 2$
	\end{case}
	
	Let $s \in S$ and $\deg(s) \geq 2$. Furthermore, let $N(s) = \{k_1, k_2, \dots, k_q\}$, where $k_i \in K$. If there exists a $j$ such that $k_j$ is not a support vertex, then $\{V_1, V_2\}$, where $V_2 = \{s, k_j\}$ and $V_1$ contains all the other vertices of $G$, forms a total transitive partition of $G$ of size 2. Therefore, we have a contradiction since $Tr_t(G) = 1$. Thus, all the vertices from $N(s) = \{k_1, k_2, \dots, k_q\}$ must be support vertices. Similarly, if $K$ has more than one vertex without a neighbour in $S$, as in Case \ref{Lower_STTCA1}, we can show a contradiction. Hence, for all $s \in S$ with $\deg(s) \geq 2$, each vertex in $N(s)$ is a support vertex, and $K$ has at most one vertex without a neighbour in $S$. This implies that $G$ satisfies condition $(b)$.
\end{proof}

Next, we characterize split graphs with $Tr_t(G)=\omega(G)-1$. Before characterizing split graphs with $Tr_t(G)=\omega(G)-1$, we have following observations for any $p$, $1\leq p\leq q-1$, where $q=\omega(G)$.

\begin{observation}\label{STTO1}
	Let $\pi=\{V_1, V_2, \ldots, V_p\}$ be a total transitive partition of a split graph $G$. Then, each $V_i$ must contain at least one vertex from $K$.
\end{observation}

\begin{observation}\label{STTO2}
	Let $\pi=\{V_1, V_2, \ldots, V_p\}$ be a total transitive partition of a split graph $G$. If $V_j$, $j\neq 1$ contains more than one vertex from $S$, then there exists $\pi'=\{V_1', V_2', \ldots, V_p'\}$ such that $|V_j'\cap S|\leq 1$ for all $j\neq 1$. Moreover, when $p=\omega(G)-1$, there exists a total transitive partition of $G$, say $\pi''=\{V_1'', V_2'', \ldots, V_p''\}$ such that $|V_i|=2$ for all $2\leq i\leq p$.
\end{observation}

\begin{observation}\label{STTO3}
	Let $\pi=\{V_1, V_2, \ldots, V_p\}$ be a total transitive partition of $G$ and $x\in S$ such that $deg(x)=1$. Then $x$ and its support vertex, say $y\in K$, cannot be in $V_t$ for $t\geq 2$.
	
\end{observation}

\begin{observation}\label{STTO4}
	Let $\pi=\{V_1, V_2, \ldots, V_p\}$ be a total transitive partition of $G$. Then there must exist a subset of $K$, say $K'$, such that $K'$ dominates $S$, and all the vertices of $K'$ must be in $V_1$. Moreover, $1\leq |K'|\leq q-p+1$, that is, the size of the $K'$ is at most $q-p+1$.
	
\end{observation}

From Observation \ref{STTO4}, we have if $\pi=\{V_1, V_2, \ldots, V_p\}$ is a total transitive partition of $G$. Then there must exist a subset of $K$, say $K'$, such that $K'$ dominates $S$, and all the vertices of $K'$ must be in $V_1$. Moreover, $1\leq |K'|\leq q-p+1$, that is, the size of the $K'$ is at most $q-p+1$. Let us define the minimum cardinality of $K'$ as \emph{S domination number by K} and denote it by $dom_K(S)$ and $K'$ with minimum cardinality as \emph{$dom_K(S)$-set}.

Let $G$ be a split graph, and $Tr_t(G)=p$ for some $1\leq p\leq \omega(G)-1$. Also, assume that $\pi=\{V_1, V_2, \ldots, V_p\}$ is a total transitive partition of $G$ with size $p$. Then we say $dom_K(S)$-set exists if there exists a subset $K'$ of $K$, such that $K'$ dominates $S$, vertices of $K'$ are in $V_1$ and $1\leq |K'|\leq q-p+1$, where $q=\omega(G)$.

In the following theorem, we characterize the split graphs with $Tr_t(G)=\omega(G)-1$.

\begin{theorem}\label{Characterization_split_upper_bound_total_transitivity}
	Let $G=(S\cup K, E)$ be a split graph, where $S$ and $K$ are the independent set and clique of $G$, respectively. Also, assume that $K$ is the maximum clique of $G$, that is, $\omega(G)=|K|$. Then $Tr_t(G)=\omega(G)-1$ if and only if $G$ satisfies one of the following conditions.
	
	\begin{enumerate}[label=(\alph*)]
		\item $dom_K(S)$-set exists and $|dom_K(S)|=1$. Also, there exist $\{s_1, \ldots, s_{\alpha}\}\subseteq S$ such that $deg(s_j)\geq j$ for all $1\leq j\leq \alpha$, where $\alpha=\omega(G)-2$ . Furthermore, for some pair of vertices $x, y$, we have an ordering of the vertices of $K\setminus \{k_1, x, y\}$, say, $\{k_2, k_3, \ldots, k_{\alpha}\}$ such that $\{k_2, k_3, \ldots, k_i\}\subseteq N(s_i)$ for all $2\leq j\leq \alpha$, where $\{k_1\}$ is a $dom_K(S)$-set.
		
		\item $dom_K(S)$-set exists and $|dom_K(S)|=2$. Also, there exist $\{s_2, \ldots, s_{\alpha}\}\subseteq S$ such that $deg(s_j)\geq j$ for all $2\leq j\leq \alpha$, where $\alpha=\omega(G)-1$ . Furthermore, we have an ordering of the vertices of $K\setminus \{k_1, k_1'\}$, say, $\{k_2, k_3, \ldots, k_{\alpha}\}$ such that $\{k_2, k_3, \ldots, k_i\}\subseteq N(s_i)$ for all $2\leq j\leq \alpha$, where $\{k_1, k_1'\}$ is a $dom_K(S)$-set.
	\end{enumerate}
\end{theorem}

\begin{proof}
	Let $G$ be a split graph, and let $G$ satisfy condition $(a)$. We show that $Tr_t(G)=\omega(G)-1$.  From the condition $(a)$, there exist $\{s_1, \ldots, s_{\alpha}\}\subseteq S$ such that $deg(s_j)\geq j$ for all $1\leq j\leq \alpha$, where $\alpha=\omega(G)-2$. Furthermore, for some pair of vertices $x, y$, we have an ordering of the vertices of $K\setminus \{k_1, x, y\}$, say $\{k_2, k_3, \ldots, k_{\alpha}\}$ such that $\{k_2, k_3, \ldots, k_i\}\subseteq N(s_i)$ for all $2\leq j\leq \alpha$, where $\{k_1\}$ is a $dom_K(S)$-set. Let us consider a vertex partition $\pi=\{V_1, V_2, \ldots, V_{\omega(G)-1}\}$ of $G$ as follows: $\{k_1, s_1\}\subseteq V_1$, $V_j=\{s_j, k_j\}$, for all $2\leq j\leq \alpha$, $V_{\omega(G)-1}=\{x, y\}$, and all the other vertices of $G$ in $V_1$.  We show that $\pi$ is a total transitive partition. Since $\{k_1\}$ is a $dom_K(S)$-set and $K$ is a clique, $V_1$ totally dominates $V_i$ for all $i\geq 1$. Again, from condition $\{k_2, k_3, \ldots, k_j\}\subseteq N(s_j)$, we have $V_s$ totally dominates $V_t$ for all $2\leq s\leq t\leq \omega(G)-2$. Furthermore, as $V_{\omega(G)-1}=\{x, y\}$ and $K$ is a clique, $V_s$ totally dominates $V_{\omega(G)-1}$ for all $1\leq s\leq \omega(G)-1$. Therefore, $\pi$ is a total transitive partition of size $\omega(G)-1$, which implies $Tr_t(G)\geq \omega(G)-1$. Moreover, from Lemma \ref{split_graph_total_transitivity_lemma_1}, we have $Tr_t(G)\leq \omega(G)-1$. Hence, if $G$ satisfies condition $(a)$, $Tr_t(G)=\omega(G)-1$.

	Now assume $G$ be a split graph and $G$ satisfies the condition $(b)$. We show that $Tr_t(G)=\omega(G)-1$.  From the condition $(b)$, there exist $\{s_2, \ldots, s_{\alpha}\}\subseteq S$ such that $deg(s_j)\geq j$ for all $2\leq j\leq \alpha$, where $\alpha=\omega(G)-1$. Also, we have an ordering of the vertices of $K\setminus \{k_1, k_1'\}$, say, $\{k_2, k_3, \ldots, k_{\alpha}\}$ such that $\{k_2, k_3, \ldots, k_i\}\subseteq N(s_i)$ for all $2\leq j\leq \alpha$, where $\{k_1, k_1'\}$ is a $dom_K(S)$-set. Let us consider a vertex partition $\pi=\{V_1, V_2, \ldots, V_{\omega(G)-1}\}$ of $G$ as follows: $\{k_1, k_1'\}\subseteq V_1$, $V_j=\{s_j, k_j\}$, for all $2\leq j\leq \alpha$, and all the other vertices of $G$ in $V_1$.  We show that $\pi$ is a total transitive partition. Since $\{k_1, k_1'\}\subseteq V_1$ is a $dom_K(S)$-set and $K$ is a clique, $V_1$ totally dominates $V_i$ for all $i\geq 1$. Again, from condition $\{k_2, k_3, \ldots, k_i\}\subseteq N(s_i)$ for all $2\leq j\leq \alpha$, we have $V_s$ totally dominates $V_t$ for all $2\leq s\leq t\leq \omega(G)-1$. Therefore, $\pi$ is a total transitive partition of size $\omega(G)-1$, which implies $Tr_t(G)\geq \omega(G)-1$. Moreover, from Lemma \ref{split_graph_total_transitivity_lemma_1}, we have $Tr_t(G)\leq \omega(G)-1$. Hence, if $G$ satisfies condition $(b)$, $Tr_t(G)=\omega(G)-1$. Therefore, if $G$ is a split graph and satisfies either $(a)$ or $(b)$ then $Tr_t(G)=\omega(G)-1$.

	Now we prove the converse part of the theorem. Let $Tr_t(G)=\omega(G)-1$. Then $dom_K(S)$-set exists and $1\leq |dom_K(S)|\leq \omega(G)-\{\omega(G)-1\}+1=2$. Based on $|dom_K(S)|=1$ or $|dom_K(S)|=2$, we divide our proof into two following cases:
	
	\begin{case}\label{Upper_STTCA1}
		$|dom_K(S)|=1$
	\end{case}
	
	By the Observation \ref{STTO2}, we can take a total transitive partition $\pi=\{V_1, V_2, \ldots, V_p\}$ of $G$ such that $|V_i|=2$ for all $2\leq i\leq p$, where $p=\omega(G)-1$. Let $\{k_1\}\subseteq K$ be a $dom_K(S)$-set. Again from Observation \ref{STTO1}, we know that each $V_i$ must contain a vertex from $K$. Therefore, exactly one set, say $V_r$ from $\pi$, contains two vertices from $K$. First, assume $V_r$ is not $V_1$, then $V_1$ contains exactly one vertex from $K$ and at least one vertex from $S$. Also, each set $V_2, \ldots, V_{r-1}, V_{r+1}, \ldots, V_p$ contains exactly one vertex from $K$ and exactly one vertex from $S$. Without loss of generality, assume $\{s_1, k_1\}\subseteq V_1$  and $V_j=\{s_j, k_j\}$ for all all $2\leq j\leq r-1, r+1\leq j\leq p$ and $V_r=\{k_r, k_r'\}$. Since $\pi$ is a total transitive partition of $G$ the $deg(s_j)\geq j$ for all $1\leq j\leq p$ and $j\neq r$. Moreover, $\{k_1, k_2, \ldots, k_j\}\subseteq N(s_j)$ for all $1\leq j\leq r-1$ and  $\{k_1, k_2, \ldots, k_{r-1}, k_{r+1}, \ldots, k_j\}\subseteq N(s_j)$ for all $r+1\leq j\leq p$. Clearly, we have the required ordering of the vertices of $K\setminus \{k_1, k_r, k_r'\}$. Similarly, we can prove that $G$ satisfies the condition $(a)$ when $V_r$ is $V_1$.
	
	\begin{case}\label{Upper_STTCA2}
		$|dom_K(S)|=2$
	\end{case}
	
	In this case, we can also assume there exists a total transitive partition $\pi=\{V_1, V_2, \ldots, V_p\}$ of $G$ such that $|V_i|=2$ for all $2\leq i\leq p$, where $p=\omega(G)-1$. Let $\{k_1, k_1'\}\subseteq K$ be a $dom_K(S)$-set. From Observation \ref{STTO1} we know that each $V_i$ must contains a vertex from $K$ and by definition of $dom_K(S)$-set $\{k_1, k_1'\}\subseteq V_1$. So, each $V_i$ contains exactly one vertex from $K$ and exactly one vertex from $S$, for all $2\leq i\leq p$. Without loss of generality assume $\{k_1, k_1'\}\subseteq V_1$  and $V_i=\{s_i, k_i\}$ for all all $1\leq i\leq p$. Since $\pi$ is a total transitive partition of $G$ the $deg(s_i)\geq j$ for all $2\leq i\leq p$. Moreover, $\{k_2, \ldots, k_i\}\subseteq N(s_i)$ for all $1\leq i\leq p$. Clearly, we have the required ordering of the vertices of $K\setminus \{k_1, k_1'\}$. 
	
	Hence, $Tr_t(G)=\omega(G)-1$ if and only if $G$ satisfies one of the above conditions.
\end{proof}

From the Lemma \ref{split_graph_total_transitivity_lemma_1}, we have $1\leq Tr_t(G)\leq \omega(G)-1$.  In the above, we already discussed the characterization of split graphs when $Tr_t(G)=1$, and $Tr_t(G)=\omega(G)-1$. Now, we show some necessary conditions for $Tr_t(G)=p$, $1\leq p\leq q-1$, where $q=\omega(G)$.

\begin{theorem}\label{Necessary_conditions_theorem_total_transitivity_split}
	Let $G=(S\cup K, E)$ be a split graph, where $S$ and $K$ are the independent set and clique of $G$, respectively. Also, assume that $K$ is the maximum clique of $G$, that is, $\omega(G)=|K|$. If $Tr_t(G)=p$, then $G$ satisfies one of the following conditions.
	
	\begin{enumerate}[label=(\alph*)]
		\item $dom_K(S)$-set exists and $|dom_K(S)|=1$
		
		\begin{enumerate}[label=(\roman*)]
			\item $|S|\geq 2p-q$
			
			\item There exists $\{s_1, s_2, \ldots, s_{\alpha}\}\subseteq S$ such that $deg(s_j)\geq j$ for all $1\leq j\leq \alpha$, where $\alpha=2p-q$
		\end{enumerate}
		
		\item $dom_K(S)$-set exists and $|dom_K(S)|\geq 2$
		
		\begin{enumerate}[label=(\roman*)]
			\item $|S|\geq 2p-q+l_{min}-2$, where $l_{min}=|dom_K(S)|$
			
			\item There exists $\{s_2, \ldots, s_{\beta}\}\subseteq S$ such that $deg(s_j)\geq j$ for all $2\leq j\leq \beta$, where $\beta=2p-q+l_{min}-2$
		\end{enumerate}
		
	\end{enumerate}
	
\end{theorem}

\begin{proof}

	Let $Tr_t(G)=p$ and $\pi=\{V_1, V_2, \ldots, V_p\}$ be a total transitive partition of $G$ with size $p$. By the Observation \ref{STTO4}, $dom_K(S)$-set exists and $1\leq |dom_K(S)|\leq p-q+1$. We divide our proof into two cases:
	
	\begin{case}\label{Any_STTCA1}
		$|dom_K(S)|=1$
	\end{case}
	
	As $\pi$ is a total transitive partition, each $|V_i|\geq 2$. Also, by Observation \ref{STTO1}, each $V_i$ must contain a vertex from $K$. Now it is given that $Tr_t(G)=p$ and $|K|=q$, so each $p-(q-p)=2p-q$ sets from $\pi$ must contain at least one vertex from $S$. Therefore, $|S|\geq p-(q-p)=2p-q$ $(a(i))$. Again, by the Observation \ref{STTO2}, each $|V_i\cap S|\leq 2$. Let $\{s_{i_1}, s_{i_2}, \ldots, s_{i_\alpha}\}\subseteq S$, where $\alpha=2p-q$ such that $s_{i_j}\in V_{i_j}$. Since $\pi$ is a total transitive partition, $deg(s_{i_j})\geq i_j\geq j$ $(a(ii))$.

	\begin{case}\label{Any_STTCA2}
		$|dom_K(S)|\geq 2$
	\end{case}
	
	In this case, let $l_{min}=|dom_K(S)|$. By the definition of $dom_K(S)$, the vertices are in $V_1$. As we know from the Observation \ref{STTO1}, each $V_i$ must contain a vertex from $K$, so each $(p-1)-[q-\{(p-1)+l_{min}\}] =2p-q+l_{min}-2$ set from $\pi$ must contain at least one vertex from $S$. Therefore, $|S|\geq 2p-q+l_{min}-2$ $(b(i))$. Again, by Observation \ref{STTO2}, each $|V_i\cap S|\leq 2$. Let $\{s_{i_1}, s_{i_2}, \ldots, s_{i_\beta}\}\subseteq S$, where $\beta=2p-q+l_{min}-2$ such that $s_{i_j}\in V_{i_j}$. Also, note that $i_j\geq 2$. Since $\pi$ is a total transitive partition, $deg(s_{i_j})\geq i_j\geq j$ $(b(ii))$.
	
	Hence, if $Tr_t(G)=p$, $G$ satisfies one of the above conditions.
\end{proof}

\begin{figure}[htbp!]
	\centering
	\includegraphics[scale=0.80]{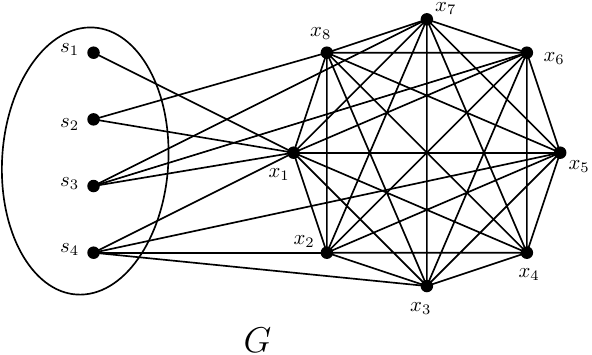}
	\caption{Example of a split graph}
	\label{fig:fsplit_graph_example_not_sufficient_condition}
\end{figure}

Next, we show that the above conditions are not sufficient for $Tr_t(G)=p$. Let us see an example of a graph given in the Figure \ref{fig:fsplit_graph_example_not_sufficient_condition}.

In the above graph $G$, $q = \omega(G) = 8$. We demonstrate that for $p = 6$, $G$ satisfies condition $(a)$ of Theorem \ref{Necessary_conditions_theorem_total_transitivity_split}, but $Tr_t(G) \neq 6$. Clearly, a $dom_K(S)$-set exists, and $dom_K(S) = \{x_1\}$, which implies $|dom_K(S)| = 1$. Also, $4 = |S| \geq 2p - q = 12 - 8$. Moreover, $\alpha = 2p - q = 4$, and $\{s_1, s_2, s_3, s_4\} \subseteq S$ such that $\deg(s_i) \geq i$. Thus, $G$ satisfies condition $(a)$ of Theorem \ref{Necessary_conditions_theorem_total_transitivity_split}.

Next, we show that $Tr_t(G) \neq 6$. Assume, for contradiction, that $Tr_t(G) = 6$, and let $\pi = \{V_1, V_2, V_3, V_4, V_5, V_6\}$ be a total transitive partition of $G$ of size $6$. By the definition of $dom_K(S)$, we know that $x_1 \in V_1$. From the degree of $s_i$, it follows that $s_i \notin V_j$ for $j \geq 5$. Since $\pi$ is a total transitive partition, $|V_i| \geq 2$, and by Observation \ref{STTO1}, each $V_i$ must contain a vertex from $K$. Therefore, each $V_i$ contains exactly one vertex from $K = \{x_1, x_2, \ldots, x_7, x_8\}$ for $1 \leq i \leq 4$. Additionally, from the degrees of $s_i$, we deduce that $s_i \in V_i$ for all $1 \leq i \leq 4$. To totally dominate $s_4$, we have $\{x_1, x_2, x_3, x_5\} \subseteq \bigcup_{i=1}^{4} V_i$, and to totally dominate $s_3$, $\{x_1, x_6, x_7\} \subseteq \bigcup_{i=1}^{3} V_i$. This leads to a contradiction, as $|K \cap \bigcup_{i=1}^{4} V_i| = 4$.

Therefore, the total transitivity of $G$ cannot be $6$. Hence, condition $(a)$ of Theorem \ref{Necessary_conditions_theorem_total_transitivity_split} is not sufficient for $Tr_t(G) = p$.

\section{NP-completeness}\label{Section_3_total}

In this section, we show that the \textsc{Maximum Total Transitivity Decision Problem} is NP-complete for a bipartite graph. Clearly, this problem is in NP. We prove the NP-hardness of this problem by showing a polynomial-time reduction from \textsc{Proper $3$-Coloring Decision Problem} in graphs to the \textsc{Maximum Total Transitivity Decision Problem}.
A \emph{proper $3$-colring} of a graph $G=(V,E)$ is a function $g: V \rightarrow \{1,2,3\}$, such that for every edge $uv \in E$, $g(u)\not= g(v) $. The proper \textsc{$3$-Coloring Decision Problem} is defined as follows:
\begin{center}

	\fbox{%
		\parbox{0.8\linewidth}{%
			\noindent\textsc{Proper $3$-Coloring Decision Problem (P$3$CDP)}
			
			\noindent\emph{Instance:} A graph $G=(V, E)$
			
			\noindent\emph{Question:} Does there exist a proper $3$-coloring of $V$?%
		}%
	}
	
\end{center}

%
%
%

The P$3$CDP is known to be NP-complete for graphs \cite{garey1990guide}. Given an instance $G=(V, E)$ of P$3$CDP, we construct an instance of \textsc{Maximum Total Transitivity Decision Problem}, that is, a graph $G'=(V', E')$ and a positive integer $k$.

\begin{figure}[!h]
	\centering
	\includegraphics[scale=0.65]{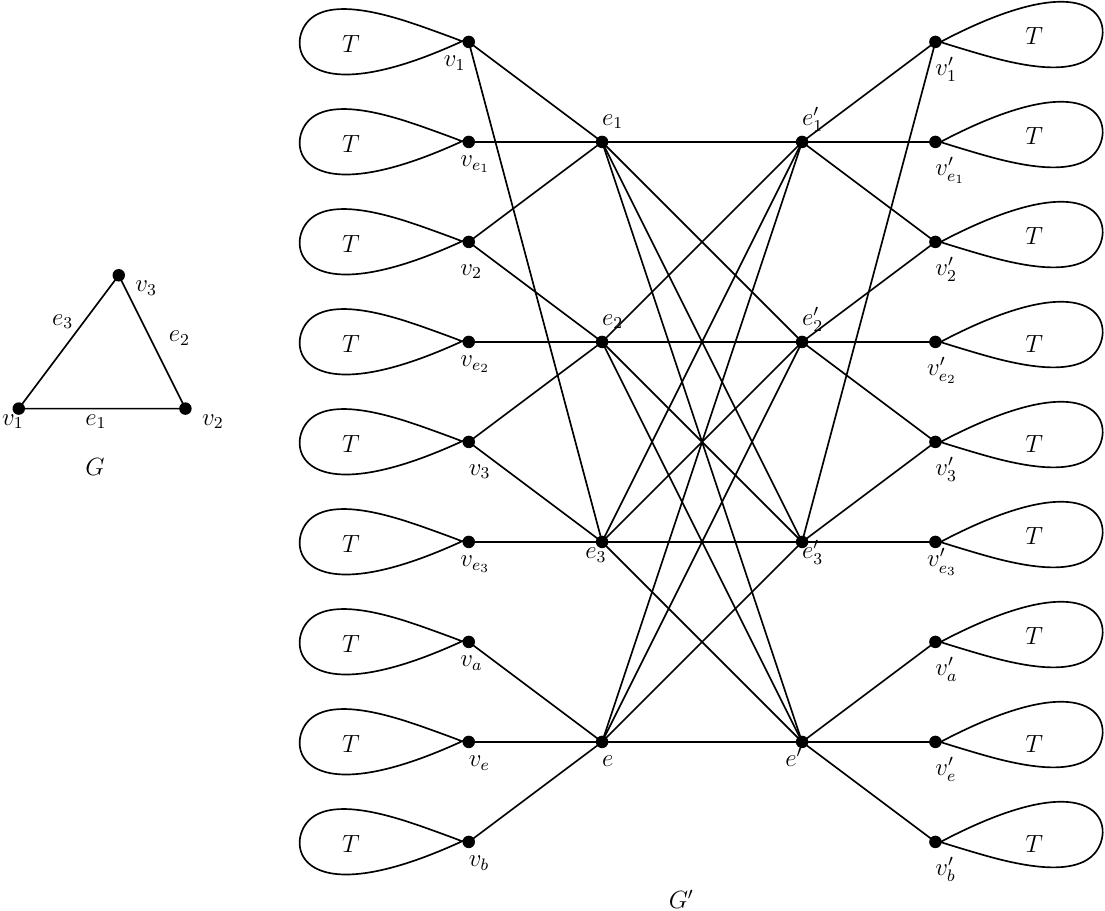}
	\caption{Construction of $G'$ from $G$ where $m=3$ and $T$ represents a $\mathcal{T}_3$}
	\label{fig:total_bipartitenp}
	
\end{figure}

The construction is as follows. Let $V=\{v_1, v_2, \ldots, v_n\}$ and $E= \{e_1, e_2, \ldots, e_m\}$. Let us consider the \emph{total transitive tree of order $3$}, that is $\mathcal{T}_3$. For each vertex $v_i\in V$, we consider two $\mathcal{T}_3$ with $v_i$ and $v_i'$ as their roots, respectively. Similarly, for each edge $e_i\in E$, we consider another two $\mathcal{T}_3$ with $v_{e_i}$ and $v_{e_i}'$ as their roots, respectively. Further, corresponding to every edge $e_j\in E$, we take two vertices $e_j, e_j'$ in $G'$. Also, we take another two vertices $e$ and $e'$ in $G'$. Let $A=\{e_1, e_2, \ldots , e_m, e\}$ and $B=\{e_1', e_2', \ldots, e_m', e'\}$. We construct a complete bipartite graph with vertex set $A\cup B$. Now, for each edge $e_t=v_iv_j\in E$, we join the edges $v_ie_t$, $v_je_t$, $v_{e_t}e_t$, $v_i'e_t'$, $v_j'e_t'$, and $v_{e_t}'e_t'$ in $G'$. Next, we consider six $\mathcal{T}_3$ with roots $v_a, v_a', v_e, v_e', v_b, v_b'$. Finally, we make $\{v_a, v_e, v_b\}$, $\{v_a', v_e', v_b'\}$ adjacent to the vertex $e$, and $e'$, respectively, and we set $k=m+4$. In this construction, we can verify that the graph $G'=(V', E')$ consists of $40m+38n+116$ vertices and $m^2+42m+34n+104$ edges. The construction is illustrated in Figure \ref{fig:total_bipartitenp}.

Next, we show that $G$ has a proper $3$-coloring if and only if $G'$ has a total transitive partition of size $k$. First, we show the forward direction of this statement in the following lemma.

\begin{lemma}
	If $G=(V,E)$ has a proper $3$-coloring, then $G'=(V', E')$ has a total transitive partition of size $k$.
\end{lemma}
\begin{proof}
	Let $g$ be a proper $3$-coloring of $G$. Based on $g$, let us consider a vertex partition, say $\pi=\{V_1,V_2,\ldots ,V_k\}$, of $G'$ as follows: for each $v_i\in V$, we put the corresponding vertex $v_i, v_i'$ in $V_q$, if $g(v_i)=q$. For each edge $e_t=v_iv_j\in E$, we put the vertex $v_{e_t}, v_{e_t}' \in V_q$, if $g(v_i)\neq q$ and $g(v_j)\neq q$. Further, we put $v_a, v_a'\in V_3$, $v_e, v_e'\in V_2$, and $v_b, v_b'\in V_1$. We put the other vertices of $\mathcal{T}_3$ according to Figure \ref{fig:tcmbt_coloring}, based on the position of the root vertex in the partition. Finally, for $1\leq j\leq m$, we put $e_{j}, e_{j}'\in V_{3+j}$ and $e, e' \in V_{m+4}$. 
	
	\begin{figure}[htbp!]
		\centering
		\includegraphics[scale=0.65]{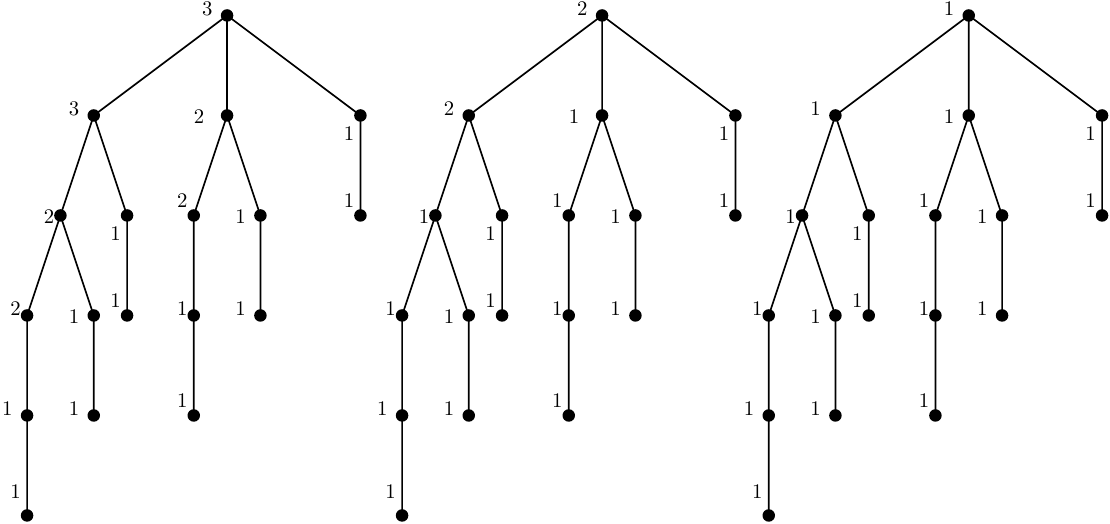}
		\caption{Partition of $\mathcal{T}_3$}
		\label{fig:tcmbt_coloring}
	\end{figure}
	
	Since the vertices of $A\cup B$ induces a complete bipartite graph, $V_i$ totally dominates $V_j$ for $4\leq i\leq j\leq k$. Also, for each $i=1,2,3$, every vertex of $A\cup B$ is adjacent to two vertices of $V_i$. Therefore, for each $i=1,2,3$, $V_i$ totally dominates $V_t$, for all $t>3$. Finally, from Figure \ref{fig:tcmbt_coloring}, it is evident that $V_i$ totally dominates $V_j$ for $1\leq i\leq j\leq 3$. Hence, $\pi$ is a total transitive partition of $G'$ of size $k$. Therefore, if $G$ has a proper $3$-coloring, then $G'$ has a total transitive partition of size $k$.
\end{proof}

Next, we show the converse of the statement is true, that is, if $G'$ has a total transitive partition of size $k$ then $G$ has a proper $3$-coloring. For this, first, we argue the existence of a total transitive partition of $G'$, say $\pi=\{V_1, V_2, \ldots, V_k\}$, such that the sets $V_1, V_2,$ and $V_3$ contain only vertices from $V'\setminus (A\cup B)$. If such a total transitive partition exists, then we can define a $3$-coloring function $g$ from $V(G)$ to $\{1, 2, 3\}$, using the partition $\pi$ as follows: $g(v)=l$, if and only if $v\in V_l$. In the following claim, we show the existence of such a total transitive partition of $G'$.

\begin{claim}\label{claim:total_transitivity_bipartite_np}
	Let  $\pi=\{V_1,V_2,\ldots ,V_k\}$ be a total transitive partition of $G'$ of size $k$ such that $|V_k|=2$. Then the sets $V_4, V_5,\ldots V_k$ contain only vertices from $A\cup B$, and the sets $V_1, V_2$, and $V_3$ contain only vertices from $V'\setminus (A\cup B)$.
\end{claim}

\begin{proof}
	First note that the vertices $\{v_a,v_e, v_b\}$, which are adjacent to $e$, must belong to $V_1, V_2$ or $V_3$. As the degree of each vertex from $\{v_a,v_e, v_b\}$ is $4$, they can be at most in $V_4$, but in this case, $e$ can be either in $V_3$ or $V_4$. Furthermore, if $e$ is in $V_3$, there is no neighbour in $V_4$ to totally dominate vertices from $\{v_a,v_e, v_b\}$ and similarly if $e$ is in $V_4$, then here is no neighbour in $V_3$ to totally  dominate vertices from $\{v_a,v_e, v_b\}$. Similarly, the vertices of $\{v'_a,v'_e,v'_b\}$, which are adjacent to $e'$ cannot be in $V_p$ for any $p\geq 4$. Further, since the degree of every vertex other than the vertices of $A\cup B$ is less than $m+4$, $V_k$ can contain two vertices only from $A\cup B=\{e_1, e_2, \ldots, e_m, e, e_1', e_2', \ldots, e_m', e'\}$. In this proof, we show that to totally dominate the vertices of $A\cup B$ which belong to $V_k$, the other vertices of $A\cup B$ must be in  $V_4, V_5,\ldots, V_{k-1}$. This, in turn, would imply that the sets $V_1, V_2$, and $V_3$ contain only vertices from $V'\setminus (A\cup B)$. The proof is divided into four cases based on whether $e\in V_k$ or not and $e'\in V_k$ or not.
	
	\begin{case}\label{TBP_case_1}
		$e\in V_k$ and $e'\in V_k$
	\end{case} 
	As discussed earlier, the vertices $\{v_a,v_e,v_b\}$ and $\{v'_a,v'_e,v'_b\}$ belong to $V_p$ for some $1\leq p \leq 3$. Since $e\in V_k$ to totally dominate $e$, each set in $\{V_1, \ldots, V_{k-1}, V_k\}$ must contain at least one vertices from $N_{G'}(e)=\{e_1', e_2', \ldots, e_m', e', v_a, v_e, v_b\}$. Since $e$ is adjacent to exactly $k=m+4$ vertices, each $V_i, 1\leq i\leq k$ contains exactly one vertex from $N_{G'}(e)$. Since $\{v_a, v_e, v_b\}$ belong to $V_p$ for some $1\leq p \leq 3$ it follows that vertices from $N_{G'}(e)\setminus \{v_a, v_e, v_b\}$ belong to $V_4, V_5,\ldots V_k$. Similarly, as $e'\in V_k$, to totally dominate $e'$, each set in $\{V_1, V_2, \ldots, V_{k-1}, V_k\}$ must contain at least one vertex from $N_{G'}(e')=\{e_1, e_2, \ldots, e, v_a', v_e', v_b'\}$. Since $e'$ is adjacent to exactly $k=m+4$ vertices, each $V_i, 1\leq i\leq k$ contains exactly one vertex from $N_{G'}(e')$. Similarly, as before, since $\{v_a', v_e', v_b'\}$ belong to $V_p$ for some $1\leq p \leq 3$, it follows that exactly three vertices from $N_{G'}(e')\setminus \{v'_a, v'_e, v'_b\}$ belong to $V_4, V_5,\ldots V_k$. Hence, the vertices of $A\cup B$ belong to $\{V_4, V_5,\ldots V_k\}$. 
	
	Clearly, none of the vertices from the set $\{v_1, \ldots, v_n, v'_1, \ldots, v'_n\}$ and the set $\{v_{e_1}, \ldots, v_{e_m}, v'_{e_1}, \ldots, v'_{e_m}\}$ can belong to $V_p$ for $p\geq 4$ because otherwise there exists a vertex in $A\cup B$ which belongs to $V_3$. Since the degree of every other vertices is at most $3$, they cannot belong to $V_p$ for $p\geq 4$. Hence, $\{V_4, V_5,\ldots V_k\}$ contain only vertices from $A\cup B$ and $\{V_1, V_2, V_3\}$ contain only vertices from $V'\setminus (A\cup B)$.
	
	\begin{case}\label{TBP_case_2}
		$e\in V_k$ and $e'\notin V_k$
	\end{case}

	Using similar arguments as in the previous case, we know that vertices from $\{v_a,v_b,v_e\}$ belong to $V_p$ for $1 \leq p \leq  3$ and the vertices from $\{e'_1, e'_2, \ldots, e'_m,e'\} $ belong to $\{V_4, V_5, \ldots, V_{k-1},V_{k}\}$ to totally dominate $e$. Moreover, every $V_i$ for $4\leq i \leq k$ contains exactly one vertex from $\{e'_1, e'_2, \ldots, e'_m, e'\}$. Since the vertices from $\{e'_1, e'_2, \ldots, e'_m, e'\}$ belong to $\{V_4, V_5, \ldots, V_{k-1},V_{k}\}$, no vertex from $\{v'_1, \ldots, v'_n, v'_{e_1}, \ldots, v'_{e_m}, v'_a, v'_b, v'_e\}$ can be in $V_p$ for $p\geq 4$. As $e'\notin V_k$, there exist a vertex from $\{e'_1, e'_2, \ldots, e'_m\} $, say $e'_r$ such that $e'_r \in V_k$. To totally dominate $e'_r$, the vertices from  $\{e_1, e_2, \ldots, e_m, e\}$  belong to $\{V_4, V_5, \ldots, V_{k}\}$. With similar arguments as in the previous case, we know that every vertex of $\{v_1, v_2, \ldots, v_n, v_{e_1}, v_{e_2}, \ldots, v_{e_m}\}$ belong to $V_p$ for $1\leq p\leq 3$. Hence, $\{V_4, V_5,\ldots V_k\}$ contain only vertices from $A\cup B$ and $\{V_1, V_2, V_3\}$ contain only vertices from $V'\setminus (A\cup B)$.

	\begin{case}\label{TBP_case_3}
		$e\notin V_k$ and $e'\in V_k$
	\end{case}
	
	Using similar arguments as in Case \ref{TBP_case_2}, we can show that $\{V_4, V_5,\ldots V_k\}$ contain only vertices from $A\cup B$ and $\{V_1, V_2, V_3\}$ contain only vertices from $V'\setminus (A\cup B)$.

	\begin{case}\label{TBP_case_4}
		$e\notin V_k$ and $e'\notin V_k$
	\end{case}
	
	Since the degree of each vertex of $V'\setminus (A\cup B)$ is at most $k-1=m+3$, they cannot belong to $V_{k}$, that is, only vertices from $A\cup B$ can be in $V_{k}$. Without loss of generality, we can assume that $e_1\in V_k$ and $e'_s \in V_{k}$ in the total transitive partition of $G'$ of size $k$. Also, let $e_1$ be the edge between $v_1$ and $v_2$ in $G$ and in the total transitive partition of $G'$ of size $k$, $v_1\in V_l$ and $v_2\in V_t$ where $l\geq t$. Next, we show that $l \leq 3$. As $v_1\in V_l$ and $\pi$ is a total transitive partition of $G$ implies that each of $V_3, V_4, \ldots, V_{l}$ contains at least one vertex from $A\setminus \{e_1, e\}$. On the other hand, as $e_s'\in V_k$, each of $\{V_{l+1}, \ldots, V_{k-1}\}$ also contains at least one vertex from $A \setminus \{e_1\}$. This implies that each of the $(m+1)$ partitions, $\{V_3,V_4, \ldots, V_{k-1}\}$, contains one vertex from the set of $m$ vertices, $A \setminus \{e_1\}$, which is impossible. Hence, $l\leq 3$. With a similar argument, we can show that $v_{e_1}$ cannot be in $V_p$ for $p\geq 4$. Therefore, each $\{V_3,V_4, \ldots, V_{k}\}$ contains exactly one vertex from $B$. Since the vertices of $B$ belong to $\{V_4, V_5, \ldots, V_{k}\}$,  the vertices of $\{v'_1, v'_2, \ldots, v'_n, v'_{e_1}, v'_{e_2}, \ldots, v'_{e_m}, v'_a, v'_b, v'_e\}$ belong to some partition $V_p$ for some $p\leq 3$. As $e'_s\in V_{k}$ and $\pi$ is a total transitive partition of $G$ to totally dominate $e'_s$, the vertices from $A\setminus \{e_1\}$ belong to $\{V_4, V_5, \ldots, V_{k-1}\}$. With similar arguments as in Case \ref{TBP_case_1}, we know that each vertex of $\{v_1, v_2, \ldots, v_n, v_{e_1}, v_{e_2}, \ldots, v_{e_m}, v_a, v_b, v_e\}$ belongs to $V_p$ for some $p\leq 3$. Hence, $\{V_4, V_5,\ldots V_k\}$ contain only vertices from $A\cup B$ and $\{V_1, V_2, V_3\}$ contain only vertices from $V'\setminus (A\cup B)$.
	
	Therefore, for all the cases, the partitions $\{V_1, V_2, V_3\}$ contain only vertices from $V'\setminus (A\cup B)$ and the partitions $\{V_4, V_5,\ldots V_k\}$ contain only vertices from $A\cup B$.
\end{proof}

Now we prove $G$ uses $3$-colors, which is a proper coloring in the following lemma.

\begin{lemma}
	If $G'$ has a total transitive partition of size $k$, then $G$ has a proper $3$-coloring.
\end{lemma}

\begin{proof}
	Let $\pi=\{V_1,V_2,\ldots,V_k\}$ be a total transitive partition of $G'$ of size $k$. By the Proposition \ref{size of total_transitive_partition} we can always assume that $|V_k|=2$.  Let us define a coloring of $G$, say $g$, by labelling $v_i$ with color $p$ if its corresponding vertex $v_i$ is in $V_p$. The previous claim ensures that $g$ is a $3$-coloring. Now we show that $g$ is a proper coloring. First note that in $G'$, since $e\in V_p$, with $p\geq 4$, $e_a', e_b'$, and $e_c'$ must belong to $V_3, V_2$, and $V_1$, respectively. Let $e_t=v_iv_j\in E$ and its corresponding vertex $e_t$ in $G'$ belong to some set $V_p$ with $p\geq 4$. This implies that the vertices $\{v_i, v_j, v_{e_t}\}$ must belong to different sets from $V_1, V_2$, and $V_3$. Therefore, $g(v_i)\neq g(v_j)$, and hence $g$ is a proper $3$-coloring of $G$.
\end{proof}

Hence, we have the following theorem.

\begin{theorem}
	The \textsc{Maximum Total Transitivity Decision Problem} is NP-complete for bipartite graphs.
\end{theorem}

\section{Bipartite chain graphs} \label{Total_tr_BCG}
In this section, we design a linear-time algorithm for finding total transitivity of a given bipartite chain graph $G$. It is known that transitivity in a bipartite chain graph $G$ is $t+1$ if and only if $G$ contains either $K_{t,t}$ or $K_{t,t}\setminus\{e\}$ as an induced subgraph for the maximum value of such $t$ \cite{paul2023transitivity}. To design the algorithm for total transitivity, we first prove that total transitivity of a bipartite chain graph $G$ can be either $Tr(G)-2$ or $Tr(G)-1$, where $Tr(G)$ is transitivity of $G$.

\begin{lemma}\label{Total_tr_BCG_lemma_1}
	Let $G=(X\cup Y,E)$ be a bipartite chain graph, and  $\sigma_X= (x_1,x_2, \ldots ,x_{n_1})$ and $\sigma_Y=(y_1,y_2, \ldots ,y_{n_2})$ be the chain ordering of $G$. Also, assume that $Tr(G)=t+1$. Then the total transitivity of $G$, $Tr(G)-2\leq Tr_t(G)\leq Tr(G)-1$.
\end{lemma}
\begin{proof}
	Since $G$ is a bipartite chain graph and $Tr(G)=t+1$, $G$ contains either $K_{t,t}$ or $K_{t,t}\setminus\{e\}$ as an induced subgraph for the maximum value of such $t$ \cite{paul2023transitivity}. Moreover, from \cite{paul2023transitivity}, we know that either $X_t\cup Y_t=V(K_{t, t})$ or $V(K_{t,t}\setminus\{e\})$ where $X_t=\{x_1, x_2, \ldots,x_t\}$ and $Y_t=\{y_1, y_2, \ldots, y_t\}$ and $\sigma_X= (x_1,x_2, \ldots ,x_m)$ and $\sigma_Y=(y_1,y_2, \ldots ,y_n)$ are the chain ordering of $G$. To prove that $Tr_t(G)\geq Tr(G)-2$, let us consider the following vertex partition $\pi=\{V_1, V_2, \ldots, V_t, V_{t-1}\}$ as follows: $\{x_1, y_1\}\subseteq V_1$, $V_j=\{x_j, y_j\}$ for all $2\leq j\leq t-1$ and all the other vertices of $G$ in $V_1$ (See Figure \ref{fig:BCG_Total_partition}).
	
	\begin{figure}[htbp!]
		\centering
		\includegraphics[scale=0.70]{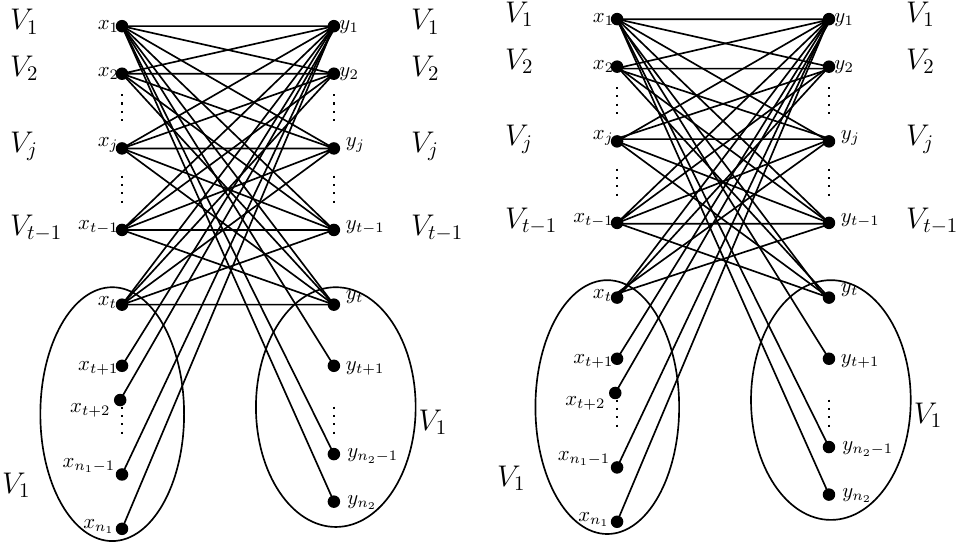}
		\caption{Left represents total transitive partition of $G$ of size at least $Tr(G)-2$ when $G$ contains $K_{t, t}$, and right represents total transitive partition of $G$ of size at least $Tr(G)-2$ when $G$ contains $K_{t, t}\setminus \{e\}$}
		\label{fig:BCG_Total_partition}
	\end{figure} 
	
	Clearly, in $\pi$, $V_i$ totally dominates $V_j$ for all $1 \leq i \leq j \leq k$. Thus, $\pi$ is a total transitive partition of $G$ with size $t-1$, that is, $Tr(G) - 2$. Therefore, $Tr_t(G) \geq Tr(G) - 2$. To prove that $Tr_t(G) \leq Tr(G) - 1$ for a bipartite chain graph, let us assume that $Tr_t(G) \geq Tr(G)$. In this case, by the definition of total transitivity, we have $Tr_t(G) = Tr(G) = t + 1$. Let $\pi = \{V_1, V_2, \ldots, V_{t+1}\}$ be a total transitive partition of $G$ of size $t + 1$. Since $|X_t| = t$, there exists at least one set in $\pi$, say $V_r$, such that $V_r$ does not contain any vertex from $X_t$. Since $\pi$ is a total transitive partition of $G$, each $V_i$ contains at least one vertex from $X$ and at least one vertex from $Y$. 
	
	Assume $x_p \in V_r$ for some $t+1 \leq p \leq n_1$. Since $V_r$ totally dominates $V_r$, there exists $y_{\alpha} \in Y_t$ such that $y_{\alpha} \in V_r$. Again, $|Y_t| = t$ and $y_{\alpha} \in Y_t \cap V_r$, implying that there exists $V_s \in \pi$ with $s \neq t$, such that $V_s$ does not contain any vertex from $Y_t$. Similarly, as $V_s$ totally dominates $V_s$, there exists $x_{\beta} \in X_t$ such that $x_{\beta} \in V_s$. Therefore, we have $\{x_{\beta}, y_q\} \subseteq V_s$ and $\{x_p, y_{\alpha}\} \subseteq V_r$. Assume the case when $s < r$. Since $\pi$ is a total transitive partition, to totally dominate $x_p$, there must be a vertex from $Y_t$ in $V_s$, which contradicts the fact that $V_s$ does not contain a vertex from $Y_t$. Similarly, we can show a contradiction when $s > r$. Therefore, our assumption that $Tr_t(G) \geq Tr(G)$ is false. Hence, for a bipartite chain graph $G$, we have $Tr(G) - 2 \leq Tr_t(G) \leq Tr(G) - 1$.
\end{proof}

Next, we characterize the bipartite chain graphs with total transitivity equal to $Tr(G)-1$.

\begin{lemma}\label{Total_tr_BCG_lemma_2}
	Let $G=(X\cup Y,E)$ be a bipartite chain graph, and  $\sigma_X= (x_1,x_2, \ldots ,x_{n_1})$ and $\sigma_Y=(y_1,y_2, \ldots ,y_{n_2})$ be the chain ordering of $G$. Also, assume that $Tr(G)=t+1$. Then $Tr_t(G)=Tr(G)-1$ if and only if $G$ contains $K_{t, t}$ as an induced subgraph for the maximum integer $t$.
\end{lemma}
\begin{proof}
	First, assume that $G$ contains $K_{t, t}$ as an induced subgraph for the maximum integer $t$. Now consider a vertex partition of $G$, say $\pi=\{V_1, V_2, \ldots, V_{t}\}$, such that $\{x_1, y_1\}\subseteq V_1$, $V_j=\{x_j, y_j\}$ for all $2\leq j\leq t$, and all the other vertices of $G$ in $V_1$. Since $G[X_t\cup Y_t]=K_{t, t}$, $V_i$ totally dominates $V_j$ for all $1\leq i\leq j\leq k$. Hence, $\pi$ is a total transitive partition of $G$. Now, by Lemma \ref{Total_tr_BCG_lemma_1}, we have $Tr_t(G)=Tr(G)-1$.

	Conversely, let $Tr_t(G)=Tr(G)-1$. As $Tr(G)=t+1$, from the results in \cite{paul2023transitivity}, we know that $G$ contains either $K_{t,t}$ or $K_{t,t}\setminus\{e\}$ as an induced subgraph for the maximum value of such $t$. If $G$ contains $K_{t,t}$ as an induced subgraph for the maximum value of such $t$, we are done. Assume the other part that is $G$ contains $K_{t,t}\setminus\{e\}$ as an induced subgraph for the maximum value of such $t$. We show that $Tr_t(G)=Tr(G)-2$, which contradicts our assumption. To show that $Tr_t(G)=Tr(G)-2$ when $G$ contains $K_{t,t}\setminus\{e\}$ as an induced subgraph for the maximum value of such $t$, we have the following claim.

	\begin{claim}
		Let $G$ be a bipartite graph such that $G$ contains $K_{t,t}\setminus\{e\}$ as an induced subgraph for the maximum value of such $t$. Then $Tr_t(G)=Tr(G)-2$
	\end{claim}
	\begin{proof}
		We show that $Tr_t(G)$ cannot be more than $Tr(G) - 2$ for a graph $G$ when $G$ contains $K_{t,t} \setminus \{e\}$ as an induced subgraph for the maximum value of such $t$. If possible, assume $Tr_t(G) \geq Tr(G) - 1$. From Lemma \ref{Total_tr_BCG_lemma_1}, we have $Tr_t(G) \leq Tr(G) - 1$. So, in this case, $Tr_t(G) = Tr(G) - 1 = t$. Let $\pi = \{V_1, V_2, \ldots, V_t\}$ be a total transitive partition of $G$. Also, assume that $x_t \in V_r$ for some $r \in \{1, 2, \ldots, t\}$. Since $\pi$ is a total transitive partition, for $x_t \in V_r$, there exists $y_{\alpha} \in \{y_1, y_2, \ldots, y_{t-1}\}$ such that $y_{\alpha} \in V_r$. Thus, there exists $V_s$, where $s \neq t$, such that $V_s$ does not contain any vertex from $\{y_1, y_2, \ldots, y_{t-1}\}$. Let $y_{\alpha'} \in V_s$ for $\alpha' \in \{t, t+1, \ldots, n_2\}$.
		
		Now, if $r < s$, then to totally dominate $y_{\alpha'}$, $V_r$ must contain a vertex from $\{x_1, x_2, \ldots, x_{t-1}\}$. Also, from the fact that $\pi$ is a total transitive partition, we have a vertex $x_{\beta} \in \{x_1, x_2, \ldots, x_{t-1}\} \cap V_r$. Therefore, in $\pi$, there exists $V_p$ for $p \neq s, r$ such that $V_p$ does not contain any vertex from $X_t$. Let $x_{\beta'} \in V_p$ with $t + 1 \leq \alpha \leq n_1$. Clearly, $p > s$, otherwise, to totally dominate $y_{\alpha'}$, there must be a vertex from $X_t$ in $V_p$. For this $x_{\beta'} \in V_p$, there must be a vertex from $\{y_1, y_2, \ldots, y_{t-1}\}$ in $V_p$. So, there exists a set $V_q$ other than $V_p$, $V_r$, and $V_s$ in $\pi$, such that it does not contain any vertices from $\{y_1, y_2, \ldots, y_{t-1}\}$. Let us assume $y_{\alpha''} \in V_q$ and $\alpha'' \in \{t, t + 1, \ldots, n_2\}$. If $p > q$, to totally dominate $x_{\beta}$, there must be a vertex from $Y_t \setminus \{y_t\}$ in $V_q$. On the other hand, if $q > p$, then to totally dominate $y_{\alpha''}$, there must be a vertex from $X_t$ in $V_p$. Both cases are not possible. Similarly, we can show contradictions when $s < r$. Therefore, our assumption is false. Hence, $Tr_t(G) = Tr(G) - 2$.
	\end{proof}
	As we assume $Tr_t(G)=Tr(G)-1$, from the above claim we can say that $G$ contains $K_{t,t}$ as an induced subgraph for the maximum value of such $t$. Hence we have the lemma.
\end{proof}

Based on Lemmas \ref{Total_tr_BCG_lemma_1} and \ref{Total_tr_BCG_lemma_2}, we can design an algorithm for finding total transitivity of a given bipartite chain graph. In \cite{paul2023transitivity}, authors showed that transitivity can be solved in linear time for a bipartite chain graph. In their algorithm they showed whether $G$ contains an induced $K_{t,t}$ or an induced $K_{t,t}\setminus\{e\}$. Hence, we have the following theorem.

\begin{theorem}
	The \textsc{Maximum Total Transitivity Problem} can be solved in linear time for bipartite chain graphs.
\end{theorem}

\section{Polynomial-time algorithm for trees}
In the Section \ref{Section_3_total}, we have seen that \textsc{Maximum Total Transitivity Decision Problem} is NP-complete for bipartite graphs. In this section, we design a polynomial-time algorithm for finding the total transitivity of a given tree $T=(V, E)$. First, we give a comprehensive description of our proposed algorithm.

\subsection{Description of the algorithm}
Let $T^c$ denote a rooted tree rooted at a vertex $c$, and $T_v^c$ denote the subtree of $T^c$ rooted at a vertex $v$. With a small abuse of notation, we use $T^c$ to denote both the rooted tree and the underlying undirected tree, whereas the total transitivity of $T^c$ is meant for the underlying undirected tree. To find the total transitivity of $T=(V, E)$, we first define the \emph{total transitive number} of a vertex $v$ in $T$. The total transitive number of a vertex $v$ in $T$ is the maximum integer $p$ such that $v\in V_p$ in a total transitive partition $\pi=\{V_1, V_2, \ldots, V_k\}$, where the maximum is taken over all total transitive partitions of $T$. We denote the total transitive number of a vertex $v$ in $T$ by $ttr(v, T)$. Note that if $Tr_t(T)=k$ and $v$ is a vertex that belongs to the set $V_k$ of a total transitive partition $\pi=\{V_1, V_2, \ldots, V_k\}$, then $ttr(v, T)$ must be equal to $k$. Therefore, the total transitivity of $T$ is the maximum total transitive number that a vertex can have; that is, $Tr_t(T)=\max\limits_{v\in V}\{ttr(v, T)\}$. Therefore, our goal is to find the total transitive number of every vertex in the tree. Next, we define another parameter, namely the \emph{rooted total transitive number}. The \emph{rooted total transitive number} of $v$ in $T^c$ is the total transitive number of $v$ in the tree $T_v^c$, and it is denoted by $ttr^r(v, T^c)$. Therefore, $ttr^r(v, T^c)= ttr(v, T_v^c)$.

To this end, we define another two parameters, namely the \emph{weak total transitive number} and \emph{rooted weak total transitive number}. The weak total transitive number of a vertex $v$ in $T$ is the maximum integer $p$ such that $v\in V_p$ in a weak total transitive partition $\pi=\{V_1, V_2, \ldots, V_k\}$. Note that the maximum is taken over all weak total transitive partitions of $T$. We denote the weak total transitive number of a vertex $v$ in $T$ by $wttr(v, T)$. Finally, we define another parameter, namely the \emph{ rooted weak total transitive number}. The \emph{rooted weak total transitive number} of $v$ in $T^c$ is the weak total transitive number of $v$ in the tree $T_v^c$, and it is denoted by $wttr^r(v, T^c)$. Therefore, $wttr^r(v, T^c)= wttr(v, T_v^c)$.

Note that the value of the rooted total transitive number of a vertex is dependent on the rooted tree, whereas the total transitive number is independent of the rooted tree. Also, for the root vertex $c$, $ttr^r(c, T^c)= ttr(c,T)$. We recursively compute the rooted total transitive number or rooted weak total transitive number of the vertices of $T^c$ in a bottom-up approach. First, we consider a vertex ordering $\sigma$, the reverse of the BFS ordering of $T^c$. For a leaf vertex $c_i$, we set $wttr^r(c_i, T^c)=1$ and for a support vertex $c_j$, $ttr^r(c_j, T^c)=wttr^r(c_j, T^c)=1$. For other non-leaf vertex $c_i$, which is not a support vertex, we call the function \textsc{Total\_Transitive\_Number$()$}, which takes the rooted total transitive number or rooted weak total transitive number of children of $c_i$ in $T^c$ as input and returns the rooted total transitive number of $c_i$ in $T^c$ or rooted weak total transitive number $c_i$ in $T^c$ according to Theorem \ref{tree_theorem_total_transitivity}. At the end of the bottom-up approach, we have the rooted total transitive number of $c$ in $T^c$ or rooted weak total transitive number of $c$ in $T^c$, that is, $ttr^r(c, T^c)$ or $wttr^r(c, T^c)$, which is the same as the total transitive number of $c$ in $T$, that is, $ttr(c, T)$ or $ttr(c, T)+1$. After the bottom-up approach, we have the total transitive number of the root vertex $c$ and the rooted total transitive number or rooted weak total transitive number of every other vertex in $T^c$. Similarly, we compute the total transitive number of all vertices of $T$ by considering a rooted tree rooted at vertex $x$. Finally, we compute the maximum overall total transitive number of vertices in $T$ and return that value as the total transitivity of $T$. The algorithm, \textsc{Total\_Transitivity(T)}, is outlined in Algorithm \ref{Algo:total_trasitivity(T)}.

\begin{algorithm}[!h]
	
	\caption{\textsc{Total\_Transitivity(T)} } \label{Algo:total_trasitivity(T)}
	
	\textbf{Input:} A tree $T=(V, E)$.
	
	\textbf{Output:} Total transitivity of $T$, that is $Tr_t(T)$.
	
	\begin{algorithmic}[1]
		
		\State Let $V=\{u_1, u_2, \ldots, u_n\}$ be the set all vertices of $T$;
		
		\ForAll {$u \in V$}
		
		\State Construct a rooted tree $T^u$, rooted at a vertex $u$;
		
		\State Let $\sigma_{u}=(c_1, c_2, \ldots, c_n=u)$ be the reverse BFS ordering of the  vertices of $T^u$, rooted at a vertex $u$;
		
		\ForAll {$c_i$ in $\sigma_{u}$}
		
		\If {$c_i$ is a leaf}
		
		\State $wttr^r(c_i, T^c)=1$;
		
		\Else
		
		\If{$c_i$ is a support vertex}
		
		\State $ttr^r(c_i, T^c)=wttr^r(c_i, T^c)=1$;
		
		\Else

		\textsc{Total\_Transitive\_Number}$(l_{i_1}', l_{i_2}', \ldots, l_{i_k}' )$.
		
		~~~~~~ /*where $l_{i_j}'=l_{i_j}$ if $l_{i_j}=ttr^r(c_{i_j}, T^u)$ otherwise $l_{i_j}'=l_{i_j}-1$*/~~~~~~
		
		~~~~~~/* where $c_{i_1}, \ldots, c_{i_k}$ are the children of $c_i$ and $l_{i_1}'\leq l_{i_2}'\leq \ldots \leq l_{i_k}'$*/~~~~
		
		\EndIf
		
		\EndIf
		
		\EndFor
		
		\If{We have rooted total transitive number of the root vertex}
		
		\State $ttr(u, T)= ttr^r(u, T^u)$;
		
		\Else
		
		\State $ttr(u, T)= wttr^r(u, T^u)-1$;
		
		\EndIf
		
		\EndFor
		
		\State $Tr_t(T)$= $\max\limits_{u_i\in V}\{ttr(u_i, T)\}$;
		
		\State \Return($Tr_t(T)$);

	\end{algorithmic}
	
\end{algorithm}

\subsection{Proof of correctness}

In this subsection, we prove the correctness of Algorithm \ref{Algo:total_trasitivity(T)}. It is clear that the correctness of Algorithm \ref{Algo:total_trasitivity(T)} depends on the correctness of the function Total\_Transitive\_Number$()$. First, we show the following two lemmas and one theorem, which prove the correctness of the Total\_Transitive\_Number$()$ function.

\begin{lemma}\label{tree_lemma_total_transitivity}
	Let $v$ be a vertex of $T$, and $ttr(v,T)=t$. Then there exists a total transitive partition of $T$, say $\{V_1, V_2, \ldots, V_i\}$, such that $v\in V_i$, for all $1\leq i\leq t$.
	
\end{lemma}

\begin{proof}
	Since $ttr(v,T)=t$, there is a total transitive partition $\pi=\{U_1, U_2, \ldots, U_t\}$ of $T$ such that $v\in U_t$. For each $1\leq i\leq t$, let us define another total transitive partition $\pi'=\{V_1,V_2,\ldots,V_i\}$ of $T$ as follows: $V_j=U_j$ for all $1\leq j\leq (i-1)$ and $V_i= \displaystyle{\bigcup_{j=i}^{t}U_j}$. Clearly, $\pi'$ is a total transitive partition of $T$ of size $i$ such that $v\in V_{i}$. Hence, the lemma follows.
\end{proof}

Similar lemma can be defined for weak total transitivity.

\begin{lemma}\label{tree_lemma_weak_total_transitivity}
	Let $v$ be a vertex of $T$, and $wttr(v,T)=t$. Then there exists a weak total transitive partition of $T$, say $\{V_1, V_2, \ldots, V_i\}$, such that $v\in V_i$, for all $1\leq i\leq t$.
	
\end{lemma}

\begin{proof}
	Since $wttr(v,T)=t$, there is a weak total transitive partition $\pi=\{U_1, U_2, \ldots, U_t\}$ of $T$ such that $v\in U_t$. Similarly, we can show that there exists a weak total transitive partition of $T$, say $\{V_1, V_2, \ldots, V_i\}$, such that $v\in V_i$, for all $1\leq i\leq t$. 
\end{proof}

In the next theorem, we find the total transitive number or weak total transitive number of a non-leaf and non-support vertex using the total transitive number or weak total transitive number of its children.

\begin{theorem}\label{tree_theorem_total_transitivity}
	Let $v_1, v_2, \ldots, v_k$ denote the children of $x$ in a rooted tree $T^x$, and for each $1\leq i\leq k$, let $l_i$ be either $ttr^r(v_i, T^x)$ or $wttr^r(v_i, T^x)$. If $l_s$ is both $ttr^r(v_s, T^x)$ and $wttr^r(v_s, T^x)$ for some $s$, we will regard $l_s$ solely as $ttr^r(v_s, T^x)$. Additionally, let $l_j' = l_j$ if $l_j = ttr^r(v_j, T^x)$; otherwise, let $l_j' = l_j - 1$, ensuring that $l_1' \leq l_2' \leq \ldots \leq l_k'$. Let $z$ denote the largest integer for which there exists a subsequence of $\{l_i': 1\leq i\leq k\}$, denoted as $(l_{i_1}'\leq l_{i_2}'\leq \ldots \leq l_{i_{z}}')$, satisfying the condition that $l_{i_{p}}'\geq p$ for all $1\leq p\leq z$. Moreover, let $q\in \{i_z+1, i_z+2, \ldots, k\}$ such that $l_q$ is a $wttr^r(v_q, T^x)$ and $l_q=l_q'+1\geq z+1$. The total transitive number or total weak transitive number of $x$ is defined in the underlying tree $T$ as follows.

	\begin{enumerate}[label=(\alph*)]
		\item If there exists $q\in \{i_z+1, i_z+2, \ldots, k\}$ such that $l_q$ is a $wttr^r(v_q, T^x)$ and $l_q=l_q'+1\geq z+1$. Then we have total transitive number of $x$ in the underlying tree, $T$, and it will be $z+1$, that is, $ttr(x, T^x)=1+z$.
		
		\item If no such $q\in \{i_z+1, i_z+2, \ldots, k\}$ exists, such that $l_q$ is a $wttr^r(v_q, T^x)$ and $l_q=l_q'+1\geq z+1$. Then we have weak total transitive number of $x$ in the underlying tree, $T$, and it will be $z+1$, that is, $wttr(x, T^x)=1+z$. 
	\end{enumerate}
\end{theorem}

\begin{proof}
	(a)First, consider the case when there exists $q \in \{i_z+1, i_z+2, \ldots, k\}$ such that $l_q$ is a $wttr^r(v_q, T^x)$ and $l_q = l_q' + 1 \geq z + 1$. Now, for each $j$, $1 \leq j \leq z$, let us consider the subtrees $T_{v_{i_j}}^x$. It is also given that $ttr^r(v_{i_j}, T^x) = l_{i_j}'$, for $j \in \{1, 2, \ldots, z\}$. For all $1 \leq p \leq z$, since $l_{i_p}' \geq p$, by Lemma \ref{tree_lemma_total_transitivity}, we know that there exists a total transitive partition $\pi^{p} = \{V_1^{p}, V_2^{p}, \ldots, V_p^{p}\}$ of $T_{v_{i_p}}^x$, such that $v_{i_p} \in V_p^{p}$. 
	
	Since $q$ exists and $wttr^r(v_q, T^x) = l_q = l_q' + 1 \geq z + 1$, by Lemma \ref{tree_lemma_weak_total_transitivity}, we know that there exists a weak total transitive partition $\pi^{q} = \{V_1^{q}, V_2^{q}, \ldots, V_{z+1}^{q}\}$ of $T_{v_q}^x$, such that $v_q \in V_{z+1}^{q}$. Let us consider a vertex partition $\pi = \{V_1, V_2, \ldots, V_z, V_{z+1}\}$ of $T^x$ as follows: $V_i = \left( \bigcup_{j=i}^{z} V_i^j \right) \cup V_i^q$, for $2 \leq i \leq z$, $V_{z+1} = \{v_q, x\}$, and all other vertices of $T^x$ are placed in $V_1$. Since $x$ is neither a leaf vertex nor a support vertex, $\pi$ is a total transitive partition of $T^x$. Therefore, by the definition of the total transitive number, $ttr(x, T^x) \geq 1 + z$.
	
	Next, we show that $ttr(x, T^x)$ cannot be more than $1 + z$. If possible, suppose $ttr(x, T^x) \geq 2 + z$. Then, by Lemma \ref{tree_lemma_total_transitivity}, there exists a total transitive partition $\pi = \{V_1, V_2, \ldots, V_{2+z}\}$ such that $x \in V_{2+z}$. This implies that for each $1 \leq i \leq 1 + z$, $V_i$ contains a neighbour of $x$, say $v_i$, such that the rooted total transitive number of $v_i$ is greater than or equal to $i$, i.e., $l_i = ttr^r(v_i, T^x) \geq i$. The set $\{l_i \mid 1 \leq i \leq 1 + z\}$ forms the desired subsequence of $\{l_i \mid 1 \leq i \leq k\}$, contradicting the maximality of $z$. Hence, $ttr(x, T^x) = 1 + z$.

	(b) Let us assume no such $q\in \{i_z+1, i_z+2, \ldots, k\}$ exits, such that $l_q$ is a $wttr^r(v_q, T^x)$ and $l_q=l_q'+1\geq z+1$. Now, for each $1\leq j\leq z$, let us consider the subtrees $T_{v_{i_j}}^x$. It is also given that $ttr^r(v_{i_j},T^x)=l_{i_j}'$, for $j\in \{1, 2, \ldots, z\}$. For all $1\leq p\leq z$, since $l_{i_{p}}'\geq p$, by Lemma \ref{tree_lemma_total_transitivity}, we know that there exists a total transitive partition $\pi^{p}=\{V_1^{p}, V_2^{p}, \ldots, V_p^{p}\}$ of $T_{v_{i_{p}}}^x$, such that $v_{i_{p}}\in V_p^{p}$. Let us consider the partition of $\pi=\{V_1, V_2, \ldots, V_z, V_{z+1}\}$ of $T^x$ as follows: $V_i=\displaystyle{\bigcup_{j=i}^{z}V_i^j}$, for $2\leq i\leq z$, $V_{z+1}=\{x\}$ and every other vertices of $T^x$ are put in $V_1$. Since $x$ is not a leaf vertex or a support vertex, $\pi$ is a weak total transitive partition of $T^x$. Therefore, $wttr(x,T)\geq 1+z$. Similarly, as above, we can show that $wttr(x, T)=1+z$.
\end{proof}

Let $x$ be neither a leaf vertex nor a support vertex of a rooted tree $T^c$. Also, assume $v_1, v_2, \ldots, v_k$ are the children of $x$ in a rooted tree $T^x$, and for each $1\leq i\leq k$, $l_i$ be either $ttr^r(v_i, T^x)$ or $wttr^r(v_i, T^x)$. If $l_s$ is both $ttr^r(v_s, T^x)$ and $wttr^r(v_s, T^x)$ for some $s$,  in that case we consider $l_s$ as  $ttr^r(v_s, T^x)$ only. Furthermore, let $l_j'=l_j$, if $l_j=ttr^r(v_j, T^x)$, otherwise, $l_j'=l_j-1$ and $l_1'\leq l_2'\leq \ldots \leq l_k'$. Now from Lemmas \ref{tree_lemma_total_transitivity}, \ref{tree_lemma_weak_total_transitivity} and Theorem \ref{tree_theorem_total_transitivity}, we have the following algorithm.

\begin{algorithm}[h]
	\caption{\textsc{Total\_Transitive\_Number$(l_1', l_2', \ldots, l_k')$} } \label{Algo:2transNumber(x)}
	
	\textbf{Input:} $l_1', l_2', \ldots, l_k'$ with $l_1'\leq l_2'\leq \ldots \leq l_k'$.
	
	\textbf{Output:} Either $ttr(x, T^x)$ or $wttr(x, T^x)$. 
	
	\begin{algorithmic}[1]
		
		\State Find the largest integer $z$ such that there exists a subsequence of $\{l_i': 1\leq i\leq k\}$, say, $(l_{i_1}'\leq l_{i_2}'\leq \ldots \leq l_{i_{z}}')$ such that  $l_{i_{p}}'\geq p$, for all $1\leq p\leq z$. 
		
		\State Find $q\in \{i_z+1, i_z+2, \ldots, k\}$ such that $l_q$ is a $wttr^r(v_q, T^x)$ and 	$l_q=l_q'+1\geq z+1$. 
		
		\If {$q$ exists}
		
		\State $ttr(x, T^x)=z+1$; 
		
		\Else
		
		\State  $wttr(x, T^x)=z+1$;
		\EndIf
		
		\State \Return($wttr(x, T^x)$ or $ttr(x, T^x)$ );

	\end{algorithmic}
	
\end{algorithm}

\subsection{Complexity Analysis}
In the function \textsc{Total\_Transitive\_Number()}, we find the total transitive number or weak total transitive number of a vertex $x$ based on the rooted total transitive number or rooted weak total transitive number of its children. We assume the children are sorted according to their $l'$-value, which is defined in Theorem \ref{tree_theorem_total_transitivity}. It is easy to see that the function \textsc{Total\_Transitive\_Number()} takes $O(deg(x))$ time for a vertex $x$. 

In the main algorithm \textsc{Total\_Transitivity(T)}, the vertex order mentioned in line $4$ can be found in linear time. Then, in a bottom-up approach, we calculate the rooted total transitive numbers or rooted weak total transitive numbers of every vertex. For that, we are spending $O(deg(c_i))$ for every $c_i\in \sigma_u$. Note that we must pass the children of $c_i$ in a sorted order to \textsc{Total\_Transitive\_Number()}. But as discussed in \cite{hedetniemi1982linear} (an algorithm for finding the Grundy number of a tree), we do not need to sort all the children based on their rooted total transitive numbers or rooted weak total transitive number; sorting the children whose rooted total transitive number or rooted weak total transitive number is less than $deg(c_i)$ is sufficient. We can argue that this can be done in $O(deg(c_i))$, as shown in \cite{hedetniemi1982linear}. Hence, the loop in lines $4-11$ takes linear time. To calculate the total transitive number for all vertices, we need $O(n^2)$ times. Therefore, we have the following theorem:

\begin{theorem}
	The \textsc{Maximum Total Transitivity Problem} can be solved in polynomial time for trees.
\end{theorem}

\section{Conclusion and future works}
In this paper, we have studied the notion of total transitivity in graphs, which is a variation of transitivity. First, we have characterized split graphs with total transitivity equal to $1$ and $\omega(G)-1$. Moreover, for split graph $G$ and $1\leq p\leq \omega(G)-1$, we have shown some necessary conditions for $Tr_t(G)=p$. Furthermore, we have shown that the decision version of this problem is NP-complete for bipartite graphs. On the positive side, we have proved that this problem can be solved in linear time for bipartite chain graphs. Finally, we have designed a polynomial-time algorithm for trees. 

This paper concludes by addressing some of the several unresolved problems in the study of total transitivity of a graph.

\begin{enumerate}
	\item Characterize split graphs with $Tr_t(G)=p$ for $p\neq 1, \omega(G)-1$.
	
	\item We know from \cite{santra2023transitivity} that in linear time we can solve the transitivity problem in split graphs. Can we design an algorithm for total transitivity in a split graph?
	
\end{enumerate}

It would be interesting to investigate the complexity status of this problem in other graph classes, such as bipartite chain graphs, bipartite permutation graphs, interval graphs, etc. Designing an approximation algorithm for this problem would be another challenging open problem.\\

\bibliographystyle{plain}
\bibliography{Total_tr_Bib}

\end{document}